\documentclass[10pt,a4paper]{amsart}
\usepackage{amssymb,amsxtra}
\usepackage[english,french]{babel}
\newtheorem{thm}{Theorem}[section]
\newtheorem{lem}[thm]{Lemma}

\newtheorem{cor}[thm]{Corollary}

\theoremstyle{definition}

\theoremstyle{remark}
\newtheorem{remark}[thm]{Remark}

\theoremstyle{plain}

\theoremstyle{remark}

\newtheorem*{example}{Example}
\numberwithin{equation}{section}

\begin{document}
\title{RESULTS ABOUT PERSYMMETRIC MATRICES OVER $ \mathbb{F}_{2}$ AND RELATED EXPONENTIALS SUMS}
\author{Jorgen~Cherly}
\address{D\'epartement de Math\'ematiques, Universit\'e de    
    Brest, 29238 Brest cedex~3, France}
\email{Jorgen.Cherly@univ-brest.fr}
\email{andersen69@wanadoo.fr}
\thanks{}
 \maketitle

  \selectlanguage{french}

  \allowdisplaybreaks
  \section{EXPONENTIAL SUMS AND RANK OF PERSYMMETRIC MATRICES OVER $ \mathbb{F}_{2}$}
  \label{sec 1}

 \begin{abstract}
 Soit  $\mathbb{K} $ le corps des s\'eries  de Laurent formelles  $ \mathbb{F}_{2}((T^{-1})). $
Nous calculons  en particulier des sommes  exponentielles dans  $\mathbb{K} $  de la forme \\
$ \sum_{degY\leq  k-1}\sum_{degZ\leq s-1}E(tYZ) $ o\`{u}
 t  est dans la boule unit\'{e} de  $\mathbb{K},$  en d\'{e}montrant qu'elles  d\'{e}pendent  seulement du rang de matrices 
 persym\'{e}triques  avec des entr\'{e}es dans $\mathbb{F}_{2} $ qui leur sont associ\'{e}es.
  ( Une matrice \;$ [\alpha _{i,j}]  $ est   persym\'{e}trique  si  $ \alpha _{i,j} = \alpha _{r,s} $ \; pour  \; i+j = r+s ). 
 En outre nous \'{e}tablissons des propri\'{e}tes de rang d'une partition de matrices persym\'{e}triques.
 Nous utilisons ces r\'{e}sultats pour calculer le nombre $\Gamma _{i}$ de matrices  persym\'{e}triques sur $\mathbb{F}_{2}$
 de rang i.
 Nous retrouvons en particulier une formule g\'{e}n\'{e}rale donn\'{e}e par D.E.Daykin.
 Notre d\' emonstration est, comme indiqu\'{e}, tr\`{e}s  diff\'{e}rente, puisqu'elle se fonde sur les propri\'{e}t\'{e}s de rang 
 d'une partition de matrices  persym\'{e}triques.
 Nous montrons \'{e}galement que le nombre R de repr\'{e}sentations dans $ \mathbb{F}_{2}[T] $ de 0 comme 
 une somme de formes quadratiques associ\'{e}es aux sommes exponentielles  $\sum_{degY\leq  k-1}\sum_{degZ\leq s-1}E(tYZ) $
  est donn\'{e} par une int\'{e}grale   \'{e}tendue \`{a}  la boule unit\'{e} et est une combinaison lin\'{e}aire des $\Gamma _{i}.$
  Nous calculons alors explicitement le nombre R.
  Des r\'{e}sultats similaires sont \'{e}galement obtenus pour les $\mathbb{K} $- espaces vectoriels de dimension n+1.
  Nous terminons notre article en calculant explicitement le nombre de matrices de rang i de
   la forme  $\left[A\over B\right], $ o\`{u} A est persym\'{e}trique.
 \end{abstract}
 
 \selectlanguage{english}

\begin{abstract}
 Let $\mathbb{K} $ be the field of Laurent Series  $ \mathbb{F}_{2}((T^{-1})). $\\
  We compute in particular   exponential sums in $\mathbb{K} $  of the form \\
   $  \sum_{degY\leq  k-1}\sum_{degZ\leq s-1}E(tYZ) $ where 
 t is in the unit interval of  $\mathbb{K},$  by showing that they only depend 
 on the rank  of some associated persymmetric matrices  with entries in  $\mathbb{F}_{2} $.
 ( A matrix\;$ [\alpha _{i,j}]  $ is  persymmetric if $ \alpha _{i,j} = \alpha _{r,s} $ \; for  \; i+j = r+s ). 
  Besides we establish rank properties of a partition of persymmetric matrices.
 We use these  results  to compute the number $ \Gamma_{i} $ of persymmetric matrices over   $ \mathbb{F}_{2} $ of
 rank i.
 We recover in this particular a general formula given by D. E. Daykin.
  Our proof is as indicated very different since it relies on rank properties 
 of a partition of persymmetric matrices. We also  prove that the number R of representations in  $ \mathbb{F}_{2}[T] $
  of 0 as a sum of some  quadratic forms associated to the exponential sums $\sum_{degY\leq  k-1}\sum_{degZ\leq s-1}E(tYZ) $
   is given by an integral over the unit interval, and is a linear combination of the $ \Gamma_{i}' s. $ We  then compute 
 explicitly the number R.
 Similar results are also obtained for n+1 dimensional $\mathbb{K} $ - vector spaces.We finish the paper by computing 
 explicitely the number of rank i matrices of the form $\left[A\over B\right], $ where A is  persymmetric.
 \end{abstract}

\subsection{An outline of the main results}
\label{subsec 1.1}

 \begin{thm}
\label{thm 1.1}
     The number $ \Gamma_{i}^{s\times k} $ of persymmetric
      $ s\times k $ matrices over $\mathbb{F}_{2}$   $$ \left ( \begin{array} {cccccc}
\alpha _{1} & \alpha _{2} & \alpha _{3} &  \ldots & \alpha _{k-1}  &  \alpha _{k} \\
\alpha _{2 } & \alpha _{3} & \alpha _{4}&  \ldots  &  \alpha _{k} &  \alpha _{k+1} \\
\vdots & \vdots & \vdots   &  \vdots  & \vdots  & \vdots \\
\vdots & \vdots & \vdots    &  \vdots & \vdots & \vdots \\
\alpha _{s-1} & \alpha _{s} & \alpha _{s+1} & \ldots  &  \alpha _{k+s-3} &  \alpha _{k+s-2}  \\
\alpha _{s} & \alpha _{s+1} & \alpha _{s+2} & \ldots  &  \alpha _{k+s-2} &  \alpha _{k+s-1}  \\
\end{array}  \right) $$ \\ of rank i  is given by  \\
 \[ \left\{\begin{array}{ccc}
             1 & if & i= 0, \\
             3\cdot 2^{2(i-1)} & if & 1\leq i \leq s-1, \\
          2^{k+s-1} - 2^{2s-2} & if &  i = s \; (s\leq k).
             \end{array}\right.\] \\
\end{thm}

\begin{remark}
\label{remark 1.2}
David E. Daykin has already proved this result over any finite field $\mathbb{F}$
with the number 2 in the formula replaced by$ |\mathbb{F}|,$   and the number 3
replaced by $ |\mathbb{F}|^2 -1 $.
Our proof is different and proper to the finite field with two elements.
\end{remark}
\begin{thm}
\label{thm 1.3}
 Let $ (j_{1},j_{2},j_{3},j_{4})\in \mathbb{N}^4 ,$  then
\begin{equation*}
{}^{\#}\Big(\begin{array}{c | c}
           j_{1} & j_{2} \\
           \hline
           j_{3} & j_{4}
           \end{array} \Big)_{\mathbb{P}/\mathbb{P}_{k+s-1}} =
  \begin{cases}
  1 & \text{if  }  j_{1}= j_{2}= j_{3}= j_{4} = 0, \\
   2^{2j-1} & \text{if  }  j_{1}=j_{2}=j_{3}=j \; ,j_{4}\in\left\{j,j+1\right\} \;,1\leq j\leq s-1,\\
    2^{2j-3} & \text {if  }  j_{1}=  j-2 \;, j_{2}= j_{3}= j-1\;,  j_{4}= j \;, 2\leq j\leq s, \\
    2^{k+s-1}-2^{2s-1} & \text{if   }  j_{1}= j_{2}= s-1 ,  j_{3}= j_{4}= s, \\
     0 & \text{otherwise},
   \end{cases}
   \end{equation*} \\
\end{thm}
  \begin{thm}
\label{thm 1.4}
 Let $ h_{s,k}(t) = h(t)  $ be the quadratic exponential sum in $ \mathbb{P} $ defined by
$$ t\in\mathbb{P}\longmapsto  \sum_{deg Y\leq k-1}\sum_{deg Z\leq s-1}E(tYZ) \in \mathbb{Z}. $$\vspace{0.5 cm}
 Then
$$ h(t) = 2^{k+s - r(D_{s\times k}(t))} $$ and \\
   $$   \int_{\mathbb{P}}h^{q}(t) dt = 2^{(q-1)(k+s) +1}\sum_{i=0}^{s}\Gamma_{i}^{s\times k}2^{-qi}.   $$
Let R denote the number of solutions
 $(Y_1,Z_1, \ldots,Y_q,Z_q) $  of the polynomial equation
                        $$ Y_1Z_1 +  Y_2Z_2 + \ldots + Y_qZ_q = 0  $$ \\
  satisfying the degree conditions \\
                   $$  degY_i \leq k-1 , \quad degZ_i \leq s-1 \quad for  1\leq i \leq q. $$ \\
 Then $$ R =  \int_{\mathbb{P}}h^{q}(t) dt $$    
     \end{thm}
\begin{thm}
\label{thm 1.5}
 Let $ g_{s,k}(t) = g(t)  $ be the quadratic exponential sum in $ \mathbb{P} $  defined by
$$ t\in\mathbb{P}\longmapsto \sum_{deg Y= k-1}\sum_{deg Z=  s-1}E(tYZ) \in \mathbb{Z}. $$
\small
  Then  \[ g(t)
 = \left\{\begin{array}{ccc}
   2^{s+k-j-2} &  if  & r(D_{(s-1) \times (k-1) }(t)) = r(D_{s \times ( k-1)}(t)) = r(D_{(s-1) \times k }(t)) = r(D_{s \times k}(t)) = j, \\
  -  2^{s+k-j-2} &  if  & r(D_{(s-1) \times (k-1) }(t)) = r(D_{s \times ( k-1)}(t)) = r(D_{(s-1) \times k }(t)) =  j \;  and \;r(D_{s \times k}(t)) = j+1, \\
   0    &   if  &  otherwise,
\end{array}\right.\] \\
\normalsize
and
\begin{equation*}
 \int_{\mathbb{P}}g^{2q}(t) dt =
      2^{(s+k-2)(2q-1)}\cdot  \sum_{j=0}^{s-1} {}^{\#}\Big(\begin{array}{c | c}
           j & j \\
           \hline
           j  &  j
           \end{array} \Big)_{\mathbb{P}/\mathbb{P}_{k +s -1}}\cdot2^{-2qj}.\\
            \end{equation*}
 \end{thm}

  \begin{thm}
 \label{thm 1.6}
 Let $ g_{m,k}(t,\eta ) = g(t,\eta ) $ be the exponential sum in $ \mathbb{P}\times\mathbb{P} $ defined by
$$ (t,\eta ) \in  \mathbb{P}\times \mathbb{P}\longmapsto 
  \sum_{deg Y\leq k-1}\sum_{deg Z\leq m}E(tYZ)\sum_{deg U =0}E(\eta YU) \in \mathbb{Z}.$$
 Then 
 \begin{equation*}
 g(t,\eta ) =  \begin{cases}
 2^{k+m+1-  r(D_{(1+m)\times k}(t)) }  & \text{if }
   r(D_{(1+m)\times k}(t)) = r(D^{\big[\stackrel{1}{1+m}\big] \times k }(t,\eta ) ), \\
     0  & \text{otherwise },
    \end{cases}
\end{equation*}
and \\
$$\int_{\mathbb{P}}\int_{\mathbb{P}}g^q(t,\eta ) dt d\eta  =
2^{q(k+m +1) -2k - m}\sum_{i = 0}^{inf(k,1+m)}
 \sigma _{i,i}^{\big[\stackrel{1}{1+m}\big]\times k}2^{-iq}. $$
 \end{thm}
  \begin{thm}
\label{thm 1.7}
 Let $ f_{m,k}(t,\eta ) = f(t,\eta ) $ be the exponential sum in $ \mathbb{P}\times\mathbb{P} $ defined by
$$ (t,\eta ) \in  \mathbb{P}\times \mathbb{P}\longmapsto 
  \sum_{deg Y\leq k-1}\sum_{deg Z\leq m}E(tYZ)\sum_{deg U \leq 0}E(\eta YU) \in \mathbb{Z}.$$
Then \\
$$ f(t,\eta ) = 2^{k+m +2 -r( r(D^{\big[\stackrel{1}{1+m}\big] \times k }(t,\eta )) }$$\\
and\\
   $$\int_{\mathbb{P}}\int_{\mathbb{P}}f^q(t,\eta ) dt d\eta  =
2^{q(k+m +2) -2k - m}\sum_{i = 0}^{inf(k,2+m)}
  \Gamma _{i}^{\Big[\substack{1 \\ 1+m }\Big] \times k}2^{-iq}. $$
 \end{thm}

    \begin{thm}
     \label{thm 1.8}
  We have the following formula  for all   $0 \leq i\leq \inf(k,2+m) $\\
  
  \begin{align*}
    \Gamma _{i}^{\Big[\substack{1 \\ 1+m }\Big] \times k}& = ( 2^{k} -  2^{i-1})\cdot\Gamma _{i-1}^{(1+m)\times k} 
 +  2^{i}\Gamma _{i}^{(1+m)\times k}. \\
  & \\
  \end{align*}

\underline {The case  k = 2}
\begin{equation*}
 \Gamma _{i}^{\Big[\substack{1 \\ 1+m }\Big] \times 2}= \begin{cases}
1 & \text{if  } i = 0, \\
 9  &  \text{if  }    i=1, \\
2^{4 +m} - 10  & \text{if   } i = 2.
\end{cases}
\end{equation*}
\underline {The case  m = 0 , $ k\geq 2 $}\\
\begin{equation*}
 \Gamma _{i}^{\Big[\substack{1 \\ 1}\Big] \times k}= \begin{cases}
 1  & \text{if  } i = 0, \\
 3\cdot(2^k-1)   &  \text{if  }    i=1, \\
2^{2k}- 3\cdot2^k + 2      & \text{if   } i = 2.
\end{cases}
\end{equation*}
\underline {The case  m = 1 , $ k\geq 3 $}\\
\begin{equation*}
 \Gamma _{i}^{\Big[\substack{1 \\ 1+ 1}\Big] \times k}= \begin{cases}
 1  & \text{if  } i = 0, \\
2^k + 5  &  \text{if  }    i=1, \\
11\cdot(2^k - 1)   & \text{if   } i = 2, \\
2^{2k+1}- 3\cdot2^{k+2} + 2^4    & \text{if   } i = 3.
\end{cases}
\end{equation*}
\underline {The case $3\leq  k \leq 1+m $}\\
\begin{equation*}
 \Gamma _{i}^{\Big[\substack{1 \\ 1+m }\Big] \times k}= \begin{cases}
1 & \text{if  } i = 0, \\
2^k +5  &  \text{if  }    i=1, \\
3\cdot2^{k+ 2i -4 } + 21\cdot2^{3i -5} & \text{if   } 2\leq i\leq k-1, \\
2^{2k +m} - 5\cdot2^{3k -5} & \text{if   } i = k.
\end{cases}
\end{equation*}
 \underline {The case $ 2\leq m\leq k-2 $}
\begin{equation*}
 \Gamma _{i}^{\Big[\substack{1 \\ 1+m }\Big] \times k}= \begin{cases}
1 & \text{if  } i = 0, \\
2^k +5  &  \text{if  }    i=1, \\
3\cdot2^{k+ 2i -4 } + 21\cdot2^{3i -5} & \text{if   } 2\leq i\leq  m, \\
11\cdot[2^{k+2m-2} - 2^{3m -2}]    & \text{if   } i = m +1,\\
2^{2k+m}  - 3\cdot2^{k+2m} +2^{3m+1}  & \text{if   } i = m +2.
\end{cases}
\end{equation*}
\end{thm}
 \begin{thm}
  \label{thm 1.9}
Let $ \Gamma _{i}^{\Big[\substack{n \\ 1+m }\Big] \times k}$ denote the number of matrices of the form
 $\left[{A\over B}\right] $of rank i such that  A is a 
$(1+m)\times k$ persymmetric matrix and B is  a $ n\times k $ matrix  over $ \mathbb{F}_{2},$ and where $\Gamma _{i}^{(1+m)\times k}$ denotes
   the number of $ (1+m)\times k $ persymmetric matrices over $ \mathbb{F}_{2}$ of rank i.\\
     Then  $\Gamma _{i}^{\Big[\substack{n \\ 1+m }\Big] \times k}$ expressed as a linear combination of the $\Gamma _{i-j}^{(1+m)\times k}$ is equal to
   \begin{equation*}
 \sum_{j= 0}^{n}2^{(n-j)\cdot(i-j)} a_{j}^{(n)}\prod_{l=1}^{j}(2^{k}- 2^{i-l})\cdot
 \Gamma _{i-j}^{(1+m)\times k} \quad for \quad 0\leq i\leq inf(k,n+m +1)
 \end{equation*}
 where 
 \begin{equation*}
 a_{j}^{(n)} = \sum_{s =0}^{j-1} (-1)^{s}\prod_{l=0}^{j-(s+1)}{2^{n+1}- 2^{l}\over 2^{j-s}-2^{l}}\cdot2^{s(n-j) +{s(s+1)\over 2}}
+ (-1)^{j}\cdot 2^{jn - {j(j-1)\over 2}}  \quad for \quad  1\leq j\leq n-1.
   \end{equation*}
   We set  
   \begin{align*}
   a_{0}^{(n)} & =   a_{n}^{(n)} = 1    \\ and\quad
   \Gamma _{i-j}^{(1+m)\times k} & = 0 \quad if \quad  i-j \notin \{0,1,2,\ldots, inf(k,1+m)\}.
   \end{align*}
  \end{thm}

\begin{cor}
\label{cor 1.10}We have the following formulas for n =  1,2,3,4,5 :\\
$ \Gamma _{i}^{\Big[\substack{1 \\ 1+m }\Big] \times k}=2^{i}\Gamma _{i}^{(1+m)\times k}+
 (2^{k}-2^{i-1})\cdot\Gamma _{i-1}^{(1+m)\times k}\quad for\quad 0\leq i\leq inf(k,2+m), $\vspace{0.5 cm}\\
  $ \Gamma _{i}^{\Big[\substack{2 \\ 1+m }\Big] \times k}=2^{2i}\Gamma _{i}^{(1+m)\times k}+
     3\cdot2^{i-1}(2^{k}-2^{i-1})\cdot\Gamma _{i-1}^{(1+m)\times k} \\
     +(2^{k}-2^{i-1})(2^{k}-2^{i-2})\cdot
      \Gamma _{i-2}^{(1+m)\times k}
       \quad for\quad 0\leq i\leq inf(k,3+m), $\vspace{0.5 cm}\\
        $ \Gamma _{i}^{\Big[\substack{3 \\ 1+m }\Big] \times k}=2^{3i}\Gamma _{i}^{(1+m)\times k}+
     7\cdot2^{(i-1)2}(2^{k}-2^{i-1})\cdot\Gamma _{i-1}^{(1+m)\times k}\\
     +   7\cdot2^{i-2}(2^{k}-2^{i-1})(2^{k}-2^{i-2})\cdot \Gamma _{i-2}^{(1+m)\times k} \\
      +(2^{k}-2^{i-1})(2^{k}-2^{i-2})(2^{k}-2^{i-3}) \Gamma _{i-3}^{(1+m)\times k}
       \quad for\quad 0\leq i\leq inf(k,4+m), $\vspace{0.5 cm}\\
        $ \Gamma _{i}^{\Big[\substack{4 \\ 1+m }\Big] \times k}=2^{4i}\Gamma _{i}^{(1+m)\times k}+
     15\cdot2^{(i-1)3}(2^{k}-2^{i-1})\cdot\Gamma _{i-1}^{(1+m)\times k}\\
     +   35\cdot2^{2i-4}(2^{k}-2^{i-1})(2^{k}-2^{i-2})\cdot\Gamma _{i-2}^{(1+m)\times k}\\
      + 15\cdot2^{i-3}(2^{k}-2^{i-1})(2^{k}-2^{i-2})(2^{k}-2^{i-3}) \Gamma _{i-3}^{(1+m)\times k}\\
       + (2^{k}-2^{i-1})(2^{k}-2^{i-2})(2^{k}-2^{i-3})(2^{k}-2^{i-4}) \Gamma _{i-4}^{(1+m)\times k}\\
        for\quad 0\leq i\leq inf(k,5+m), $\vspace{0.5 cm}\\
       $ \Gamma _{i}^{\Big[\substack{5 \\ 1+m }\Big] \times k}=2^{5i}\Gamma _{i}^{(1+m)\times k}+
     31\cdot2^{(i-1)4}(2^{k}-2^{i-1})\cdot\Gamma _{i-1}^{(1+m)\times k}\\
     +   155\cdot2^{3i-6}(2^{k}-2^{i-1})(2^{k}-2^{i-2})\cdot\Gamma _{i-2}^{(1+m)\times k}\\
      + 155\cdot2^{2i-6}(2^{k}-2^{i-1})(2^{k}-2^{i-2})(2^{k}-2^{i-3}) \Gamma _{i-3}^{(1+m)\times k}\\
       +31\cdot2^{i-4}(2^{k}-2^{i-1})(2^{k}-2^{i-2})(2^{k}-2^{i-3})(2^{k}-2^{i-4}) \Gamma _{i-4}^{(1+m)\times k} \\
       + (2^{k}-2^{i-1})(2^{k}-2^{i-2})(2^{k}-2^{i-3})(2^{k}-2^{i-4})(2^{k}-2^{i-5}) \Gamma _{i-5}^{(1+m)\times k}\\
        for\quad 0\leq i\leq inf(k,6+m). $\vspace{0.5 cm}\\
    \end{cor}  
   
\begin{thm}
\label{thm 1.11}
 Let $ \displaystyle  f_{m,k}(t,\eta_{1},\eta _{2},\ldots,\eta _{n} ) $  be the exponential sum  in $ \mathbb{P}^{n+1} $ defined by\\
  $ \displaystyle (t,\eta_{1},\eta _{2},\ldots,\eta _{n} )\in \mathbb{P}^{n+1}\longrightarrow \\
  \sum_{deg Y\leq k-1}\sum_{deg Z\leq m}E(tYZ)\sum_{deg U_{1}\leq  0}E(\eta_{1} YU_{1})
  \sum_{deg U_{2} \leq 0}E(\eta _{2} YU_{2}) \ldots \sum_{deg U_{n} \leq 0} E(\eta _{n} YU_{n}). $\vspace{0.5 cm}\\
  Set $$(t,\eta_{1},\eta _{2},\ldots,\eta _{n} ) =
  \big(\sum_{i\geq 1}\alpha _{i}T^{-i}, \sum_{i\geq 1}\beta  _{1i}T^{-i},\ldots, \sum_{i\geq 1}\beta  _{ni}T^{-i}) \in\mathbb{P}^{n+1}.   $$           
   Then
  $$ f_{m,k}(t,\eta_{1},\eta _{2},\ldots,\eta _{n} ) = 
  2^{k+m +n+1-r( D^{\Big[\stackrel{n}{1+m}\Big] \times k }(t,\eta_{1},\eta _{2},\ldots,\eta _{n} ))} $$
 where 
$$   D^{\Big[\stackrel{n}{1+m}\Big] \times k }(t,\eta_{1},\eta _{2},\ldots,\eta _{n} ) $$
 denotes  the following  $(1+n+m)\times k $ matrix
   $$  \left ( \begin{array} {cccccc}
\alpha _{1} & \alpha _{2} & \alpha _{3} &  \ldots & \alpha _{k-1}  &  \alpha _{k} \\
\alpha _{2 } & \alpha _{3} & \alpha _{4}&  \ldots  &  \alpha _{k} &  \alpha _{k+1} \\
\vdots & \vdots & \vdots   &  \ldots  & \vdots  &  \vdots \\
\vdots & \vdots & \vdots    &  \ldots & \vdots  &  \vdots \\
\alpha _{1+m} & \alpha _{2+m} & \alpha _{3+m} & \ldots  &  \alpha _{k+m-1} &  \alpha _{k+m}  \\
\hline
\beta  _{11} & \beta  _{12} & \beta  _{13} & \ldots  &  \beta_{1 k-1} &  \beta _{ 1k}  \\
\beta  _{21} & \beta  _{22} & \beta  _{23} & \ldots  &  \beta_{2 k-1} &  \beta _{2k}\\
\vdots & \vdots & \vdots   &  \ldots  & \vdots  &  \vdots \\
\vdots & \vdots & \vdots    &  \ldots & \vdots  &  \vdots \\
\beta  _{n1} & \beta  _{n2} & \beta  _{n3} & \ldots  &  \beta_{n k-1} &  \beta _{nk}
 \end{array}  \right). $$\vspace{0.5 cm}
Then the number denoted by $ R_{q}(n,k,m) $ of solutions \\
 $(Y_1,Z_1,U_{1}^{(1)},U_{2}^{(1)}, \ldots,U_{n}^{(1)}, Y_2,Z_2,U_{1}^{(2)},U_{2}^{(2)}, 
\ldots,U_{n}^{(2)},\ldots  Y_q,Z_q,U_{1}^{(q)},U_{2}^{(q)}, \ldots,U_{n}^{(q)}   ) $ \vspace{0.5 cm}\\
 of the polynomial equations  \vspace{0.5 cm}
  \[\left\{\begin{array}{c}
 Y_{1}Z_{1} +Y_{2}Z_{2}+ \ldots + Y_{q}Z_{q} = 0  \\
   Y_{1}U_{1}^{(1)} + Y_{2}U_{1}^{(2)} + \ldots  + Y_{q}U_{1}^{(q)} = 0  \\
    Y_{1}U_{2}^{(1)} + Y_{2}U_{2}^{(2)} + \ldots  + Y_{q}U_{2}^{(q)} = 0\\
    \vdots \\
   Y_{1}U_{n}^{(1)} + Y_{2}U_{n}^{(2)} + \ldots  + Y_{q}U_{n}^{(q)} = 0 
 \end{array}\right.\]
    satisfying the degree conditions \\
                   $$  degY_i \leq k-1 , \quad degZ_i \leq m ,
                   \quad degU_{j}^{i} \leq 0 , \quad  for \quad 1\leq j\leq n  \quad 1\leq i \leq q $$ \\
  is equal to the following integral over the unit interval in $ \mathbb{K}^{n+1} $
 $$ \int_{\mathbb{P}^{n+1}} f_{m,k}^{q}(t,\eta_{1},\eta _{2},\ldots,\eta _{n} )dt d\eta_{1}d\eta _{2}\ldots d\eta _{n}. $$
  Observing that $ f_{m,k}(t,\eta_{1},\eta _{2},\ldots,\eta _{n} )$ is constant on cosets of $ \mathbb{P}_{k+m}\times\mathbb{P}_{k}^{n}, $\;
  the above integral is equal to 
$$  2^{q(k+m +n+1) -(n+1)k - m}\sum_{i = 0}^{inf(k,n+1+m)}
   \Gamma _{i}^{\left[\stackrel{n}{1+m}\right]\times k}2^{-iq} = R_{q}(n,k,m)$$
\end{thm}
\begin{example}
\label{example 1.12}
The number $R_{q}(0,k,m)$  of solutions
 $(Y_1,Z_1, \ldots,Y_q,Z_q) $  of the polynomial equation
                        $$ Y_1Z_1 +  Y_2Z_2 + \ldots + Y_qZ_q = 0  $$ \\
  satisfying the degree conditions \\
                   $$  degY_i \leq k-1 , \quad degZ_i \leq m \leq k-1 \quad for  1\leq i \leq q. $$ \\
 is equal to the following integral   \\ 
  \begin{align*}
& \int_{\mathbb{P}}\Big[ \sum_{deg Y\leq k-1}\sum_{deg Z\leq m}E(tYZ) \Big]^{q}dt
  = 2^{(q-1)(k+m+1) +1}\sum_{i=0}^{1+m}\Gamma_{i}^{(1+m)\times k}2^{-qi}  \\
 &  = \left\{\begin{array}{ccc}
   2^k + 2^{1+m} -1 & if & q=1,\\
   2^{2k} + 3\cdot(m+1)\cdot 2^{k+m} & if & q=2, \\
   2^{(q-1)(k+m+1) +1}\left[1+ 3\frac{1-2^{(2-q)m}}{2^{q}- 2^{2}}+
 (2^{k+m}-2^{2m})2^{-q(1+m)}\right]
 & if & 3\leq q.
  \end{array}\right. 
  \end{align*}
\end{example}

\begin{example}
\label{example 1.13}
The number  $ \Gamma _{i}^{\Big[\substack{1 \\ 1+2 }\Big] \times 3}$ of rank i matrices of the form \\
  $$  \left ( \begin{array} {ccc}
\alpha _{1} & \alpha _{2} & \alpha _{3}  \\
\alpha _{2 } & \alpha _{3} & \alpha _{4} \\
\alpha _{3} & \alpha _{4} & \alpha _{5} \\
\hline
\beta  _{1} & \beta  _{2} & \beta  _{3} \\
 \end{array}  \right) $$
is equal to 

\begin{equation*}
  \begin{cases}
1 & \text{if  } i = 0, \\
13 &  \text{if  }    i=1, \\
66  & \text{if   } i = 2, \\
176  & \text{if   } i = 3.
\end{cases}
\end{equation*}
\vspace{0.1 cm}

 The number $R_{q}(1,3,2)$ of solutions 
 $(Y_1,Z_1,U_{1}, \ldots,Y_q,Z_q,U_{q}) $  of the polynomial equations
   \[\left\{\begin{array}{c}
 Y_{1}Z_{1} +Y_{2}Z_{2}+ \ldots + Y_{q}Z_{q} = 0  \\
   Y_{1}U_{1} + Y_{2}U_{2} + \ldots  + Y_{q}U_{q} = 0
 \end{array}\right.\]
  satisfying the degree conditions \\
                   $$  degY_i \leq 2 , \quad degZ_i \leq 2 ,\quad degU_{i}\leq 0 \quad for \quad 1\leq i \leq q $$ 
  is equal to the following integral \\
   \begin{align*}
  & \int_{\mathbb{P}}\int_{\mathbb{P}}\big[  \sum_{deg Y\leq 2}\sum_{deg Z\leq 2}E(tYZ)\sum_{deg U \leq 0}E(\eta YU) \big]^{q}dtd\eta  
  =  2^{7q-8}\sum_{i = 0}^{3}\Gamma _{i}^{\Big[\substack{1 \\ 1+2 }\Big] \times 3}2^{-iq}\\
  & \\
  & = 2^{4q-8}\cdot\big[2^{3q} +13\cdot2^{2q} +66\cdot2^{q} +176 \big].
   \end{align*}
\end{example}

\begin{example}
\label{example 1.14}
 The number  $ \Gamma _{i}^{\Big[\substack{5 \\ 1+2 }\Big] \times 4}$ of rank i matrices of the form \\
  $$  \left ( \begin{array} {cccc}
\alpha _{1} & \alpha _{2} & \alpha _{3} & \alpha _{4} \\
\alpha _{2 } & \alpha _{3} & \alpha _{4} & \alpha _{5}\\
\alpha _{3} & \alpha _{4} & \alpha _{5} & \alpha _{6} \\
\hline
\beta  _{11} & \beta  _{12} & \beta  _{13} & \beta_{14}  \\
\beta  _{21} & \beta  _{22} & \beta  _{23} & \beta_{24}  \\
\beta  _{31} & \beta  _{32} & \beta  _{33} & \beta_{34}  \\
\beta  _{41} & \beta  _{42} & \beta  _{43} & \beta_{44}  \\
\beta  _{51} & \beta  _{52} & \beta  _{53} & \beta_{54} 
 \end{array}  \right) $$
is equal to 

\begin{equation*}
  \begin{cases}
1 & \text{if  } i = 0, \\
561 &  \text{if  }    i=1, \\
65670 & \text{if   } i = 2, \\
3731208 & \text{if   } i = 3, \\
63311424   & \text{if   } i = 4. 
\end{cases}
\end{equation*}
\vspace{0.1 cm}

 The number $R_{3}(5,4,2)$ of solutions \\
 
  $(Y_1,Z_1,U_{1}^{(1)},U_{2}^{(1)},U_{3}^{(1)},U_{4}^{(1)},U_{5}^{(1)}, Y_2,Z_2,U_{1}^{(2)},U_{2}^{(2)}, 
U_{3}^{(2)},U_{4}^{(2)},U_{5}^{(2)}, Y_3,Z_3,U_{1}^{(3)},U_{2}^{(3)},U_{3}^{(3)},U_{4}^{(3)},U_{5}^{(3)}   ) $ \\

 of the polynomial equations
  \[\left\{\begin{array}{c}
 Y_{1}Z_{1} +Y_{2}Z_{2} + Y_{3}Z_{3} = 0,  \\
   Y_{1}U_{1}^{(1)} + Y_{2}U_{1}^{(2)}  + Y_{3}U_{1}^{(3)} = 0,  \\
 Y_{1}U_{2}^{(1)} + Y_{2}U_{2}^{(2)}  + Y_{3}U_{2}^{(3)} = 0, \\
 Y_{1}U_{3}^{(1)} + Y_{2}U_{3}^{(2)}  + Y_{3}U_{3}^{(3)} = 0, \\
 Y_{1}U_{4}^{(1)} + Y_{2}U_{4}^{(2)}  + Y_{3}U_{4}^{(3)} = 0, \\
  Y_{1}U_{5}^{(1)} + Y_{2}U_{5}^{(2)}  + Y_{3}U_{5}^{(3)} = 0, 
 \end{array}\right.\]
   satisfying the degree conditions \\
                   $$  degY_i \leq 3 , \quad degZ_i \leq 2 ,\quad degU_{i}\leq 0 \quad  for \quad 1\leq j\leq 5  \quad 1\leq i \leq 3 $$ \\ 
    is equal  to the following integral over the unit interval in $ \mathbb{K}^{6} $\\
   \begin{align*}
  & \int_{\mathbb{P}^{6}} f_{2,4}^{3}(t,\eta_{1},\eta _{2},\eta _{3},\eta _{4},\eta _{5} )dt d\eta_{1}d\eta _{2}d\eta _{3}d\eta _{4}d\eta _{5}
  =  2^{10}\cdot \sum_{i = 0}^{4}\Gamma _{i}^{\Big[\substack{5 \\ 1+2 }\Big] \times 4}2^{-i3}= 24413824.
   \end{align*}
   
    \end{example}
    \newpage
    
 \section{ EXPONENTIAL SUMS AND RANK OF DOUBLE  PERSYMMETRIC MATRICES OVER $\mathbb{F}_{2} $}
 \label{sec 2}
   \selectlanguage{french}
\begin{abstract}
 Soit $\mathbb{K}^{2} $  le  $\mathbb{K}$ - espace vectoriel de dimension 2  o\`{u} $\mathbb{K}$  d\'{e}note 
    le corps des s\'eries  de Laurent formelles  $ \mathbb{F}_{2}((T^{-1})). $ Nous calculons en particulier des sommes exponentielles 
    (dans $\mathbb{K}^{2}$) de la forme \\
     $  \sum_{\deg Y \leq  k-1}\sum_{\deg Z \leq s-1}E(tYZ)\sum_{\deg U \leq s+m-1}E(\eta YU) $
   o\`{u} $ (t,\eta ) $ est dans la boule  unit\'{e} de $\mathbb{K}^{2}.$\\
     Nous d\'{e}montrons qu'elles d\'{e}pendent uniquement
   du rang de matrices doubles persym\'{e}triques  avec des entr\'{e}es dans $\mathbb{F}_{2},$
    c'est-\`{a}-dire des matrices de la forme  $\left[A\over B\right]  $ o\`{u}  A  est une matrice  $ s \times k $
    persym\'{e}trique et B une matrice  $ (s+m) \times k $  persym\'{e}trique  (une matrice \;$ [\alpha _{i,j}]  $ est 
      persym\'{e}trique  si  $ \alpha _{i,j} = \alpha _{r,s} $ \; pour  \; i+j = r+s). En outre, nous \'{e}tablissons plusieurs 
      formules concernant des propri\'{e}t\'{e}s de rang de partitions de matrices doubles persym\'{e}triques, ce qui nous 
      conduit \`{a} une formule r\'{e}currente du nombre   $  \Gamma _{i}^{\Big[\substack{s \\ s+m }\Big] \times k} $ 
   des matrices de rang i de la forme  $\left[A\over B\right]. $ Nous d\'{e}duisons de cette formule r\'{e}currente 
   que si  $  0\leq i\leq\inf(s-1, k-1),$ le nombre    $  \Gamma _{i}^{\Big[\substack{s \\ s+m }\Big] \times k} $ d\'{e}pend uniquement de i.
   D'autre part, si   $ i\geq s+1, k\geq i,\; \Gamma _{i}^{\Big[\substack{s \\ s+m }\Big] \times k}  $ peut \^{e}tre  calcul\'{e}
   \`{a} partir du nombre  $  \Gamma _{s' +1}^{\Big[\substack{s' \\ s'+m' }\Big] \times k'} $de matrices  de rang (s'+1) de la forme  
       $\left[A'\over B'\right] $ o\`{u}  A' est une matrice   $ s' \times k' $ persym\'{e}trique  et  B' une matrice  $ (s'+m') \times k' $ persym\'{e}trique,
        o\`{u}  s', m' et k' d\'{e}pendent  de  i, s, m et k. La preuve de ce r\'{e}sultat est bas\'{e}e sur une formule (donn\'{e}e dans  [4])
        du nombre de matrices de rang i de la forme     $\left[A\over b_{-}\right] $   o\`{u}  A est  persym\'{e}trique  et  $ b_{-} $
 une matrice ligne avec entr\'{e}es dans $ \mathbb{F}_{2}.$ Nous montrons \'{e}galement que le nombre R de repr\'{e}sentations dans 
  $ \mathbb{F}_{2}[T] $ 
 des \'{e}quations polynomiales
   \[\left\{\begin{array}{cc}
YZ + Y_{1}Z_{1} + \ldots + Y_{q-1}Z_{q-1}= 0 \\
YU + Y_{1}U_{1} + \ldots + Y_{q-1}U_{q-1}= 0
    \end{array}\right.\]\\ 
  associ\'{e}es aux sommes exponentielles \\
   $\sum_{degY\leq  k-1}\sum_{degZ\leq s-1}E(tYZ)\sum_{\deg U\leq s+m-1}E(\eta YU)  $
  est donn\'{e} par une int\'{e}grale sur la boule unit\'{e} de $\mathbb{K}^{2}$ et est une combinaison lin\'{e}aire de
  $ \Gamma _{i}^{\Big[\substack{s \\ s+m }\Big] \times k}$ pour $ i\geq 0. $ Nous pouvons alors calculer explicitement le nombre R.
 \end{abstract}
\allowdisplaybreaks
 \selectlanguage{english}

 \begin{abstract}
Let $ \mathbb{K}^{2}  $ be the 2-dimensional vectorspace over $ \mathbb{K}$  where $ \mathbb{K}$ denotes the field of
 Laurent Series  $ \mathbb{F}_{2}((T^{-1})). $ We compute in particular   exponential sums, (in $\mathbb{K}^{2} $)  of the form \\
  $  \sum_{\deg Y \leq  k-1}\sum_{\deg Z \leq s-1}E(tYZ)\sum_{\deg U \leq s+m-1}E(\eta YU) $ where 
$ (t,\eta ) $ is in the unit interval of  $\mathbb{K}^{2}.$  We  show that they only depend  on the rank  of some associated  double  persymmetric matrices  
with entries in  $\mathbb{F}_{2}, $ that is matrices of the form $\left[A\over B\right] $ 
 where A is  a  $ s \times k $ persymmetric matrix and B a  $ (s+m) \times k $ persymmetric matrix. 
  (A matrix $\; [\alpha _{i,j}]  $ is  persymmetric if $ \alpha _{i,j} = \alpha _{r,s}  \; for  \; i+j = r+s ).$
   Besides, we establish several formulas concerning  rank properties of  partitions of  double  persymmetric matrices, which leads to a 
    recurrent formula for the number $  \Gamma _{i}^{\Big[\substack{s \\ s+m }\Big] \times k} $ of rank i  matrices of the form $ \left[{A\over B}\right].$  
    We deduce from the recurrent formula that  if  $  0\leq i\leq\inf(s-1, k-1) $ then   $  \Gamma _{i}^{\Big[\substack{s \\ s+m }\Big] \times k} $ 
     depends only on  i. On the other hand,  if $ i\geq s+1, k\geq i,\; \Gamma _{i}^{\Big[\substack{s \\ s+m }\Big] \times k}  $
      can be computed from the number  $  \Gamma _{s' +1}^{\Big[\substack{s' \\ s'+m' }\Big] \times k'} $ of rank (s'+1) matrices of the form 
       $\left[A'\over B'\right] $ 
 where A' is  a  $ s' \times k' $ persymmetric matrix and B' a  $ (s'+m') \times k' $ persymmetric matrix, where s', m' and k' depend on i, s, m and k.
 The proof of this result is based on a formula (given in [4]) of the number of rank i matrices of the form  $\left[A\over b_{-}\right] $ 
 where A is persymmetric and $ b_{-} $
 a one-row matrix with entries in  $ \mathbb{F}_{2}.$
   We also  prove that the number R of representations in  $ \mathbb{F}_{2}[T] $ 
 of the polynomial equations
   \[\left\{\begin{array}{cc}
YZ + Y_{1}Z_{1} + \ldots + Y_{q-1}Z_{q-1}= 0 \\
YU + Y_{1}U_{1} + \ldots + Y_{q-1}U_{q-1}= 0
    \end{array}\right.\]\\ 
  associated to the exponential sums \\
   $\sum_{degY\leq  k-1}\sum_{degZ\leq s-1}E(tYZ)\sum_{\deg U\leq s+m-1}E(\eta YU)  $
   is given by an integral over the unit interval of $\mathbb{K}^{2}$, and is a linear combination of the
    $ \Gamma _{i}^{\Big[\substack{s \\ s+m }\Big] \times k}\;for\;  i\geq 0 . $ We can  then compute 
 explicitly the number R.
   \end{abstract}

 \subsection{An outline of the main results}
\label{subsec 2.1}

\begin{thm}
\label{thm 2.1}Let q be a rational integer $ \geq 1,$ then \\
 \begin{align}
    g_{k,s,m}(t,\eta )  =  g(t,\eta ) & = \sum_{deg Y\leq k-1}\sum_{deg Z \leq  s-1}E(tYZ)\sum_{deg U \leq s+m-1}E(\eta YU)  =
   2^{2s+m+k- r( D^{\left[\stackrel{s}{s+m}\right] \times k }(t,\eta  ) )}, \label{eq 2.1} \\
   & \nonumber \\
    \int_{\mathbb{P}\times \mathbb{P}} g_{k,s,m}^{q}(t,\eta )dtd\eta 
 & =2^{(2s+m+k)(q-1)}\cdot 2^{-k+2}\cdot \sum_{i = 0}^{\inf(2s+m,k)} \Gamma _{i}^{\Big[\substack{s  \\ s+m }\Big] \times k}\cdot2^{- qi}. \label{eq 2.2} \\
  & \nonumber 
 \end{align}
\end{thm}

  \begin{thm}
 \label{thm 2.2} Let $ s\geq 2, \; m\geq 0, \; k\geq 1 $ and  $ 0\leq i\leq \inf{(2s+m,k)}. $ Then we have the following recurrent formula for the number
  $ \Gamma _{i}^{\Big[\substack{s \\ s+m }\Big] \times k} $ of rank i matrices 
 of the form $ \left[{A\over B}\right],$ $$   \left ( \begin{array} {cccccc}
\alpha _{1} & \alpha _{2} & \alpha _{3} &  \ldots & \alpha _{k-1}  &  \alpha _{k} \\
\alpha _{2 } & \alpha _{3} & \alpha _{4}&  \ldots  &  \alpha _{k} &  \alpha _{k+1} \\
\vdots & \vdots & \vdots    &  \vdots & \vdots  &  \vdots \\
\alpha _{s-1} & \alpha _{s} & \alpha _{s +1} & \ldots  &  \alpha _{s+k-3} &  \alpha _{s+k-2}  \\
\alpha _{s} & \alpha _{s+1} & \alpha _{s +2} & \ldots  &  \alpha _{s+k-2} &  \alpha _{s+k-1}  \\
\hline \\
\beta  _{1} & \beta  _{2} & \beta  _{3} & \ldots  &  \beta_{k-1} &  \beta _{k}  \\
\beta  _{2} & \beta  _{3} & \beta  _{4} & \ldots  &  \beta_{k} &  \beta _{k+1}  \\
\vdots & \vdots & \vdots    &  \vdots & \vdots  &  \vdots \\
\beta  _{m+1} & \beta  _{m+2} & \beta  _{m+3} & \ldots  &  \beta_{k+m-1} &  \beta _{k+m}  \\
\vdots & \vdots & \vdots    &  \vdots & \vdots  &  \vdots \\
\beta  _{s+m-1} & \beta  _{s+m} & \beta  _{s+m+1} & \ldots  &  \beta_{s+m+k-3} &  \beta _{s+m+k-2}  \\
\beta  _{s+m} & \beta  _{s+m+1} & \beta  _{s+m+2} & \ldots  &  \beta_{s+m+k-2} &  \beta _{s+m+k-1} 
\end{array}  \right). $$  such that A is a $ s\times k $ persymmetric matrix
  and B a  $ (s+m)\times k $ persymmetric matrix  with entries in  $ \mathbb{F}_{2}: $\\
   \begin{align}
  \Gamma _{i}^{\Big[\substack{s \\ s+m }\Big] \times k} & = 2\cdot \Gamma _{i-1}^{\Big[\substack{s -1\\ s-1+(m+1) }\Big] \times k}
+ 4\cdot \Gamma _{i-1}^{\Big[\substack{s \\ s+(m-1) }\Big] \times k} - 8\cdot \Gamma _{i-2}^{\Big[\substack{s-1 \\ s-1+m }\Big] \times k}
 + \Delta _{i}^{\Big[\substack{s \\ s+m }\Big] \times k} \label{eq 2.3} \\
 & \text{where the remainder  $\Delta _{i}^{\Big[\substack{s \\ s+m }\Big] \times k}$ is equal to} \nonumber \\
  & \sigma _{i,i,i}^{\left[\stackrel{s-1}{\stackrel{s+m-1 }{\overline {\stackrel{\alpha_{s -}}{\beta_{s+m-} }}}}\right] \times k }
  - 3\cdot \sigma _{i-1,i-1,i-1}^{\left[\stackrel{s-1}{\stackrel{s+m-1 }{\overline {\stackrel{\alpha_{s -}}{\beta_{s+m-} }}}}\right] \times k } 
  + 2\cdot \sigma _{i-2,i-2,i-2}^{\left[\stackrel{s-1}{\stackrel{s+m-1 }{\overline {\stackrel{\alpha_{s -}}{\beta_{s+m-} }}}}\right] \times k }.\label{eq 2.4} 
  \end{align}
    Recall that   $$  \sigma _{i,i,i}^{\left[\stackrel{s-1}{\stackrel{s+m-1 }
{\overline {\stackrel{\alpha_{s -}}{\beta_{s+m-} }}}}\right] \times k } $$\\  is equal to  the cardinality of the following set \\
 \small
 $$   \left\{(t,\eta )\in \mathbb{P}/\mathbb{P}_{k+s-1}\times \mathbb{P}/\mathbb{P}_{k+s+m-1}
\mid  r(D^{\big[\stackrel{s-1}{s-1+m}\big] \times k }(t,\eta ))  = r(D^{\big[\stackrel{s}{s+m-1}\big] \times k }(t,\eta ))
= r(D^{\big[\stackrel{s}{s+m}\big] \times k }(t,\eta )) = i
 \right\}. $$\\ 
   \end{thm}
\begin{thm}
\label{thm 2.3}Let $ s\geq 2 $ and $ m\geq 0, $ we have in the following two cases : \vspace{0.1 cm}\\
\underline {The case $1 \leq  k \leq 2s+m-2 $}\\
\begin{equation}
\label{eq 2.5}
 \sigma _{i,i,i}^{\left[\stackrel{s}{\stackrel{s+m }{\overline {\stackrel{\alpha_{s -}}{\beta_{s+m-} }}}}\right] \times k }= 
  \begin{cases}
1 & \text{if  } i = 0,\quad k\geq 1, \\
   4 \cdot \Gamma _{i}^{\Big[\substack{s -1\\ s-1 +m }\Big] \times i} -  \Gamma _{i+1}^{\Big[\substack{s-1 \\ s-1+m }\Big] \times (i+1)} 
      & \text{if   } 1 \leq i \leq k-1, \\
   4 \cdot \Gamma _{k}^{\Big[\substack{s -1\\ s-1 +m }\Big] \times k}    & \text{if   } i = k.
\end{cases}
\end{equation}

\underline {The case $ k\geq  2s+m-2 $ }\\
\begin{equation}
\label{eq 2.6}
 \sigma _{i,i,i}^{\left[\stackrel{s}{\stackrel{s+m }{\overline {\stackrel{\alpha_{s -}}{\beta_{s+m-} }}}}\right] \times k }= 
  \begin{cases}
1 & \text{if  } i = 0,          \\
   4 \cdot \Gamma _{i}^{\Big[\substack{s -1\\ s-1 +m }\Big] \times i} -  \Gamma _{i+1}^{\Big[\substack{s-1 \\ s-1+m }\Big] \times (i+1)} 
      & \text{if   } 1 \leq i \leq 2s+m-3,  \\
   4 \cdot \Gamma _{2s+m-2}^{\Big[\substack{s -1\\ s-1 +m }\Big] \times (2s+m-2)}    & \text{if   } i = 2s+m-2. 
\end{cases}
\end{equation}
\end{thm}

\begin{thm}The remainder $\Delta _{i}^{\Big[\substack{s \\ s+m }\Big] \times k} $in the recurrent formula is equal to \vspace{0.1 cm}\\
\label{thm 2.4}
   \begin{equation}
\label{eq 2.7}
  \begin{cases}
1 & \text{if  } i = 0,\; k \geq 1,         \\
   4 \cdot \Gamma _{1}^{\Big[\substack{s -1\\ s-1 +m }\Big] \times 1} -  \Gamma _{2}^{\Big[\substack{s-1 \\ s-1+m }\Big] \times 2} 
      & \text{if   } i = 1,\; k\geq 2, \\
    4 \cdot \Gamma _{1}^{\Big[\substack{s -1\\ s-1 +m }\Big] \times 1} - 3   & \text{if   } i = 1,\; k =1, \\
     7 \cdot \Gamma _{2}^{\Big[\substack{s -1\\ s-1 +m }\Big] \times 2} -  12\cdot \Gamma _{1}^{\Big[\substack{s-1 \\ s-1+m }\Big] \times 1}
     -  \Gamma _{3}^{\Big[\substack{s -1\\ s-1 +m }\Big] \times 3} + 2  & \text{if   } i = 2,\; k\geq 3, \\
        7 \cdot \Gamma _{2}^{\Big[\substack{s -1\\ s-1 +m }\Big] \times 2} -  12\cdot \Gamma _{1}^{\Big[\substack{s-1 \\ s-1+m }\Big] \times 1} +2
         & \text{if   } i = 2,\; k = 2, \\
         7 \cdot \Gamma _{i}^{\Big[\substack{s -1\\ s-1 +m }\Big] \times i} -  14\cdot \Gamma _{i-1}^{\Big[\substack{s-1 \\ s-1+m }\Big] \times (i-1)}
     +8\cdot \Gamma _{i-2}^{\Big[\substack{s -1\\ s-1 +m }\Big] \times (i-2)}  -  \Gamma _{i+1}^{\Big[\substack{s -1\\ s-1 +m }\Big] \times (i+1)}   & \text{if } 3\leq i\leq 2s+m-3,  k\geq i+1, \\
        7 \cdot \Gamma _{i}^{\Big[\substack{s -1\\ s-1 +m }\Big] \times i} -  14\cdot \Gamma _{i-1}^{\Big[\substack{s-1 \\ s-1+m }\Big] \times (i-1)}
     +8\cdot \Gamma _{i-2}^{\Big[\substack{s -1\\ s-1 +m }\Big] \times (i-2)}    & \text{if   }  3\leq i\leq 2s+m-3, \; k = i, \\
       7 \cdot \Gamma _{2s+m-2}^{\Big[\substack{s -1\\ s-1 +m }\Big] \times (2s+m-2)} -  14\cdot \Gamma _{2s+m-3}^{\Big[\substack{s-1 \\ s-1+m }\Big] \times (2s+m-3)}
     +8\cdot \Gamma _{2s+m-4}^{\Big[\substack{s -1\\ s-1 +m }\Big] \times (2s+m-4)}    & \text{if   }  i = 2s+m-2, \; k \geq  i, \\
     - 14\cdot \Gamma _{2s+m-2}^{\Big[\substack{s-1 \\ s-1+m }\Big] \times (2s+m-2)}
     +8\cdot \Gamma _{2s+m-3}^{\Big[\substack{s -1\\ s-1 +m }\Big] \times (2s+m-3)}    & \text{if   }  i = 2s+m-1, \; k \geq  i, \\
      8\cdot \Gamma _{2s+m-2}^{\Big[\substack{s -1\\ s-1 +m }\Big] \times (2s+m-2)}    & \text{if   }  i = 2s+m , \; k \geq  i. 
\end{cases}
\end{equation}
  \end{thm}

     \begin{thm}
\label{thm 2.5}
We have
  \begin{align}
\Delta _{i}^{\Big[\substack{s \\ s+m }\Big] \times k}&  =  \Delta _{i}^{\Big[\substack{s \\ s+m }\Big] \times (i+1)} & \;  for \; i\in [0,2s+m-3],  \; k\geq i +1, \label{eq 2.8}\\
 \Delta _{i}^{\Big[\substack{s \\ s+m }\Big] \times k} & =  \Delta _{i}^{\Big[\substack{s \\ s+m }\Big] \times i} & \;  for \; i\in \left\{2s+m-2, 2s+m-1,2s+m \right\}, \; k\geq i. \label{eq 2.9}
\end{align}
\end{thm}

  \begin{thm}
 \label{thm 2.6}We have for all $ m\geq 0 $\vspace{0.01 cm}\\
  \begin{equation}
\label{eq 2.10}
 \Gamma _{j}^{\Big[\substack{s \\ s +m }\Big] \times (k+1)}
  - \Gamma _{j}^{\Big[\substack{s \\ s +m}\Big] \times k}= 0  \quad \text{if $ \quad 0\leq j\leq s-1,\;k >j. $}
 \end{equation}
  We have in the cases  $ m \in \left\{0,1\right\} $ \vspace{0.01 cm}\\
 \begin{equation}
\label{eq 2.11}
 \Gamma _{s+j}^{\Big[\substack{s \\ s}\Big] \times (k+1)}
  - \Gamma _{s+j}^{\Big[\substack{s \\ s}\Big] \times k}=
  \begin{cases}
 3\cdot 2^{k+s-1}  & \text{if  } j = 0,\; k >s, \\
 21\cdot 2^{k+s+3j-4}   &  \text{if  }    1\leq j\leq s-1,\; k>s+j, \\
3\cdot2^{2k+2s-2}- 3\cdot2^{k+4s - 4}       & \text{if   } j  = s,  \; k>2s,
\end{cases}
\end{equation}
 \begin{equation}
\label{eq 2.12}
 \Gamma _{s+j}^{\Big[\substack{s \\ s+1}\Big] \times (k+1)}
  - \Gamma _{s+j}^{\Big[\substack{s \\ s+1}\Big] \times k}=
  \begin{cases}
  2^{k+s-1}  & \text{if  } j = 0,\; k >s, \\
  11\cdot2^{k+s-1}& \text{if  } j = 1,\; k > s +1, \\
21 \cdot 2^{k+s+3j-5}   &  \text{if  }    2\leq j\leq s ,\; k>s+j, \\
3\cdot2^{2k+2s-1}- 3\cdot2^{k+4s -2}       & \text{if   } j  = s+1 ,  \; k>2s +1.
\end{cases}
\end{equation}
 In the case $ m\geq 2  $\vspace{0.1 cm}\\
 \begin{equation}
\label{eq 2.13}
 \Gamma _{s+j}^{\Big[\substack{s \\ s+m}\Big] \times (k+1)}
  - \Gamma _{s+j}^{\Big[\substack{s \\ s+m}\Big] \times k}=
  \begin{cases}
  2^{k+s-1}  & \text{if  } j = 0,\; k >s, \\
  3\cdot2^{k+s +2j -3}& \text{if  } 1\leq j\leq m-1,\; k > s+j,  \\
11\cdot2^{k+ s + 2m -3}              &  \text{if  }  j = m, \; k > s+m,
\end{cases}
\end{equation}
\begin{equation}
\label{eq 2.14}
 \Gamma _{s+m+j}^{\Big[\substack{s \\ s+m}\Big] \times (k+1)}
  - \Gamma _{s+m+j}^{\Big[\substack{s \\ s+m}\Big] \times k}=
  \begin{cases}
  21\cdot2^{k+s +2m +3j -4}& \text{if  } 1\leq j\leq s-1,\; k > s+m+ j,  \\
3\cdot2^{2k+2s+m-2} -  3\cdot2^{k+4s+2m-4}            &  \text{if  }  j = s, \; k > 2s+m.
\end{cases}
\end{equation}
 \end{thm}
 
 \begin{thm}
 \label{thm 2.7}We have for $ m\geq 1 $ \vspace{0.01 cm}\\
   \begin{align}
  \Gamma _{s +j}^{\Big[\substack{s \\ s +m }\Big] \times k} &  =  8^{j-1}\cdot \Gamma _{s+1}^{\Big[\substack{s \\ s +(m-(j-1)) }\Big] \times (k-(j-1))}  && \text{if\;$1\leq j\leq m,\;k\geq s+j, $} \label{eq 2.15}\\
  & \nonumber \\
    \Gamma _{s+1}^{\Big[\substack{s \\ s +(m-(j-1)) }\Big] \times (s+1)} & =  2^{4s+(m-(j-1))} - 3\cdot2^{3s -1} + 2^{2s -1} &&\text{if\; $ 1\leq j\leq m,\; k = s+j, $} \label{eq 2.16}\\ 
    & \nonumber \\
     \Gamma _{s+1}^{\Big[\substack{s \\ s +(m-(j-1)) }\Big] \times (k-(j-1))}& = 3\cdot2^{k - j +s} +21\cdot[2^{3s-1}-2^{2s-1}] && \text{if\;$1\leq j\leq m-1,\;k>s+j, $} \label{eq 2.17} \\
      & \nonumber \\
  \Gamma _{s +1}^{\Big[\substack{s \\ s + 1 }\Big] \times (k-(m-1))}& = 11\cdot2^{k- m +s} + 21\cdot2^{3s- 1} - 11\cdot2^{2s- 1} && \text{if\;$ j = m,\;k>s+m. $} \label{eq 2.18} \\
   & \nonumber 
\end{align}
 \end{thm}
  \begin{thm}
\label{thm  2.8}We have for $ m\geq 0 $ \vspace{0.1 cm}
\begin{align}
 \Gamma _{s +m +1 +j}^{\Big[\substack{s \\ s +m }\Big] \times k} &  =  8^{2j + m}\cdot \Gamma _{s-j+1}^{\Big[\substack{s -j\\ s -j }\Big] \times (k- m-2j)}  &&\text{if   \; $ 0\leq j\leq s-1,\;  k\geq s + m+1+j $}, \label{eq 2.19}\vspace{0.1 cm} \\
  & \nonumber \\
 \Gamma _{s - j+1}^{\Big[\substack{s -j\\ s -j }\Big] \times ( s - j+1)} & = 2^{4s-4j} -3\cdot2^{3s-3j-1} + 2^{2s-2j-1} &&\text{if   \; $ 0\leq j\leq s-1,\;  k =  s + m+1+j $}, \label{eq 2.20} \vspace{0.1 cm} \\
  & \nonumber \\
 \Gamma _{s-j+1}^{\Big[\substack{s -j\\ s -j }\Big] \times (k- m-2j)} 
   & = 21\cdot[2^{k-m -3j+s-1} + 2^{3s-3j-1} -5\cdot2^{2s-2j-1}] &&\text{ if $ \; 0\leq j\leq s-2,\; k > s+m+ 1+j $}, \label{eq 2.21}\vspace{0.1 cm} \\
    & \nonumber \\
    \Gamma _{2}^{\Big[\substack{1 \\ 1  }\Big] \times (k-m -2s +2)}  & =  2^{2(k-m) -4s +4} -3\cdot2^{k-m -2s +2} +2  &&\text{if   \; $ j =  s-1,\;  k > 2s + m $}. \label{eq 2.22} \\
     & \nonumber 
  \end{align}
\end{thm}  
 \begin{thm}
\label{thm  2.9}  We have \vspace{0.1 cm} \\
  \begin{equation}
\label{eq 2.23} \Gamma _{i}^{\Big[\substack{s \\ s  }\Big] \times k} = 
\begin{cases}
1 & \text{if  } i = 0,\; k \geq 1,         \\
 21\cdot2^{3i-4}  - 3\cdot2^{2i-3}  & \text{if   } 1\leq i\leq s-1,\; k > i, \\
 3\cdot2^{k+s-1} + 21\cdot2^{3s-4} - 27\cdot2^{2s-3}  & \text{if   }  i = s ,\; k > s, \\
 21\cdot[ 2^{k - 2s +3i - 4}  + 2^{3i -4} - 5\cdot2^{4i-2s -5} ] & \text{if   } s+1\leq i\leq 2s-1,\; k > i, \\
 2^{2k + 2s  -2}  -3\cdot 2^{k+4s - 4}  + 2^{6s  -5} & \text{if   }  i = 2s ,\; k > 2s. 
  \end{cases}    
   \end{equation}
     \end{thm} 
     
 \begin{thm}
\label{thm  2.10}  We have \vspace{0.1 cm} \\
  \begin{equation}
\label{eq 2.24} \Gamma _{i}^{\Big[\substack{s \\ s  }\Big] \times i} = 
\begin{cases}
2^{2s +2i -2}  - 3\cdot2^{3i-4} + 2^{2i-3} & \text{if   } 1\leq i\leq s,  \\
2^{2s +2i -2}  - 3\cdot2^{3i-4} + 2^{4i -2s -5} & \text{if   } s+1\leq i\leq 2s. 
  \end{cases}    
   \end{equation} 
   \end{thm} 
    
 \begin{thm}
\label{thm 2.11}  We have \vspace{0.1 cm} \\
  \begin{equation}
\label{eq 2.25} \Gamma _{i}^{\Big[\substack{s \\ s +1 }\Big] \times k} = 
\begin{cases}
1 & \text{if  } i = 0,\; k \geq 1,         \\
 21\cdot2^{3i-4}  - 3\cdot2^{2i-3}  & \text{if   } 1\leq i\leq s-1,\; k > i, \\
 2^{k+s-1} + 21\cdot2^{3s-4} - 11\cdot2^{2s-3}  & \text{if   }  i = s ,\; k > s, \\
 11\cdot 2^{k+s-1} + 21\cdot2^{3s- 1} - 53 \cdot2^{2s - 1}  & \text{if   }  i = s +1 ,\; k > s +1, \\
 21\cdot[ 2^{k - 2s +3i -  5}  + 2^{3i -4} - 5\cdot2^{4i-2s - 6} ] & \text{if   } s+ 2\leq i\leq 2s ,\; k > i, \\
 2^{2k + 2s  - 1}  -3\cdot 2^{k+4s - 2}  + 2^{6s  - 2} & \text{if   }  i = 2s +1 ,\; k > 2s +1.
  \end{cases}    
   \end{equation}
     \end{thm} 
     
 \begin{thm}
\label{thm  2.12}  We have  \vspace{0.1 cm} \\
  \begin{equation}
\label{eq 2.26} \Gamma _{i}^{\Big[\substack{s \\ s +1 }\Big] \times i} = 
\begin{cases}
2^{2s +2i - 1}  - 3\cdot2^{3i-4} + 2^{2i-3} & \text{if   } 1\leq i\leq s +1, \\
2^{2s +2i -2}  - 3\cdot2^{3i-4} + 2^{4i -2s - 6} & \text{if   } s+2\leq i\leq 2 s +1.
  \end{cases}    
   \end{equation} 
   \end{thm} 
    
 \begin{thm}
\label{thm  2.13}  We have for $ m\geq 2 $ \vspace{0.1 cm} \\
  \begin{equation}
\label{eq 2.27} \Gamma _{i}^{\Big[\substack{s \\ s + m }\Big] \times k} = 
\begin{cases}
1 & \text{if  } i = 0,\; k \geq 1,         \\
 21\cdot2^{3i-4}  - 3\cdot2^{2i-3}  & \text{if   } 1\leq i\leq s-1,\; k > i, \\
 2^{k+s-1} + 21\cdot2^{3s-4} - 11\cdot2^{2s-3}  & \text{if   }  i = s ,\; k > s, \\
 3\cdot2^{k-s+2i-3} + 21\cdot[2^{3i-4} - 2^{3i -s-4}] & \text{if   }  s+1 \leq i\leq s+m-1 ,\; k > i,\\
 11\cdot 2^{k+s+2m-3} + 21\cdot2^{3s +3m-4} - 53 \cdot2^{2s  +3m- 4}  & \text{if   }  i = s +m ,\; k > s +m, \\
 21\cdot[ 2^{k - 2s +3i - m -4}  + 2^{3i -4} - 5\cdot2^{4i-2s - m-5} ] & \text{if   } s +m +1\leq i\leq 2s +m -1 ,\; k > i, \\
 2^{2k + 2s  +m -2}  -3\cdot 2^{k+4s +2m-4}  + 2^{6s  + 3m -5} & \text{if   }  i = 2s + m ,\; k > 2s + m.
  \end{cases}    
   \end{equation}
     \end{thm} 
     
 \begin{thm}
\label{thm  2.14}  We have for $ m\geq 2  $ \vspace{0.1 cm} \\
  \begin{equation}
\label{eq 2.28} \Gamma _{i}^{\Big[\substack{s \\ s +m }\Big] \times i} = 
\begin{cases}
2^{2s +2i +m -2}  - 3\cdot2^{3i-4} + 2^{2i-3} & \text{if   } 1\leq i\leq s +1, \\
2^{2s +2i  +m -2 }  - 3\cdot2^{3i-4} + 2^{3i - s - 4} & \text{if   } s+2\leq i\leq s+m+1,\\
2^{2s +2i  +m -2 }  - 3\cdot2^{3i-4} + 2^{4i - 2s - m -5 } & \text{if   } s+m+2\leq i\leq 2s+m. 
  \end{cases}    
   \end{equation} 
   \end{thm}

 \begin{thm}
\label{thm 2.15}
 We denote by  $ R_{q}(k,s,m) $ the number of solutions \\
 $(Y_1,Z_1,U_{1}, \ldots,Y_q,Z_q,U_{q}) $  of the polynomial equations
   \[\left\{\begin{array}{c}
 Y_{1}Z_{1} +Y_{2}Z_{2}+ \ldots + Y_{q}Z_{q} = 0,  \\
   Y_{1}U_{1} + Y_{2}U_{2} + \ldots  + Y_{q}U_{q} = 0,
 \end{array}\right.\]
  satisfying the degree conditions \\
                   $$  degY_i \leq k-1 , \quad degZ_i \leq s-1 ,\quad degU_{i}\leq s+m-1 \quad for \quad 1\leq i \leq q. $$ \\                           
Then \\
\begin{align}
  R_{q}(q,k,s,m) & =  \int_{\mathbb{P}\times \mathbb{P}} g_{k,s,m}^{q}(t,\eta )dtd\eta  \label{eq 2.29}\\
&   = 2^{(2s+m+k)(q-1)}\cdot 2^{-k+2}\cdot \sum_{i = 0}^{\inf(2s+m,k)} \Gamma _{i}^{\Big[\substack{s  \\ s+m }\Big] \times k}\cdot2^{- qi}. \nonumber 
\end{align}
  \end{thm}
\begin{example} Let q = 3, k = 4, s = 3, m = 2. Then 
   \begin{equation*}
   \Gamma _{i}^{\Big[\substack{3 \\  3 + 2 }\Big] \times 4} = 
\begin{cases}
1 & \text{if  } i = 0,        \\
 9  & \text{if   } i = 1, \\
 78   & \text{if   }  i = 2, \\
 648   & \text{if   }  i = 3,   \\
 15648 &  \text{if   }  i = 4. 
  \end{cases}    
   \end{equation*}\vspace{0.1 cm}\\
Hence  the number $ R_{3}(4,3,2) $ of solutions \\
 $(Y_1,Z_1,U_{1}, Y_2,Z_2,U_{2},Y_3,Z_3,U_{3}) $  of the polynomial equations
   \[\left\{\begin{array}{c}
 Y_{1}Z_{1} +Y_{2}Z_{2} + Y_{3}Z_{3} = 0,  \\
   Y_{1}U_{1} + Y_{2}U_{2}  + Y_{3}U_{3} = 0,
 \end{array}\right.\]
  satisfying the degree conditions \\
                   $$  degY_i \leq 3 , \quad degZ_i \leq 2 ,\quad degU_{i}\leq 4 \quad for \quad 1\leq i \leq 3 $$ \\                           
  is equal to \vspace{0.1 cm}\\
      \begin{align*}
   &  \int_{\mathbb{P}\times \mathbb{P}} g_{4,3,2}^{3}(t,\eta )dtd\eta 
 =  2^{22}\cdot \sum_{i = 0}^{4} \Gamma _{i}^{\Big[\substack{3 \\ 3 +2 }\Big] \times 4}\cdot2^{-3i}  \\
 & =  2^{22}\cdot[1 + 9\cdot2^{-3} + 78\cdot2^{-6} + 648\cdot2^{-9} + 15648\cdot2^{-12}] \\
&  = 35356672.
\end{align*}
\end{example}
\begin{example} Let  q = 4, k = 6, s = 5, m = 0. Then  
   \begin{equation*}
   \Gamma _{i}^{\Big[\substack{5 \\  5}\Big] \times 6} = 
\begin{cases}
1 & \text{if  } i = 0,        \\
 9  & \text{if   } i = 1, \\
 78   & \text{if   }  i = 2, \\
 648   & \text{if   }  i = 3,   \\
 5280 &  \text{if  }  i = 4, \\
 42624 & \text{if  } i  = 5, \\
  999936  &  \text{if   }  i = 6.
  \end{cases}    
   \end{equation*}\vspace{0.1 cm}\\
   Hence  the number  $  R_{4}(6,5,0) $ of solutions  \vspace{0.1 cm}\\
 $(Y_1,Z_1,U_{1}, Y_2,Z_2,U_{2},Y_3,Z_3,U_{3}, Y_4, Z_4,U_{4} ) $  of the polynomial equations
   \[\left\{\begin{array}{c}
 Y_{1}Z_{1} +Y_{2}Z_{2} + Y_{3}Z_{3} + Y_{4}Z_{4} = 0,  \\
   Y_{1}U_{1} + Y_{2}U_{2}  + Y_{3}U_{3} + Y_{4}U_{4}  = 0,
 \end{array}\right.\]
  satisfying the degree conditions \\
                   $$  degY_i \leq 5 , \quad degZ_i \leq 4 ,\quad degU_{i}\leq 4 \quad for \quad 1\leq i \leq 4 $$ \\                           
  is equal to \vspace{0.1 cm}\\
     \begin{align*}
&   \int_{\mathbb{P}\times \mathbb{P}} g_{4,5,0}^{4}(t,\eta )dtd\eta \\
&  =  2^{44}\cdot \sum_{i = 0}^{6} \Gamma _{i}^{\Big[\substack{5 \\  5 }\Big] \times 6}\cdot2^{-4i}  \\
 & =  2^{44}\cdot[1 + 9\cdot2^{-4} + 78\cdot2^{-8} + 648\cdot2^{-12} + 5280 \cdot2^{-16} + 42624\cdot 2^{-20} + 999936\cdot 2^{-24}] \\
&  = 37014016\cdot 2^{20}.
\end{align*}
\end{example}  
 \begin{example} The fraction of square  double persymmetric  $\left[s\atop s+m\right]\times (2s+m) $ matrices which are 
 invertible is equal to $ \displaystyle \frac{ \Gamma _{2s+m}^{\left[s\atop s+m\right] \times (2s+m)}}{\sum_{i = 0}^{2s+m}\Gamma _{i}^{\left[s\atop s+m\right] \times (2s+m)} } = \frac{3}{8}. $
   \end{example}
   
  \newpage

\section{Exponential sums and rank of triple persymmetric matrices over $\mathbb{F}_{2}$ }
\label{sec 3}
 \begin{abstract}
  \selectlanguage{french} 
 
 Notre travail concerne  une g\' en\' eralisation des r\' esultats obtenus dans :
 {Exponential sums and rank of  double  persymmetric  matrices over  $\mathbf{F}_2 $  }\\
{arXiv : 0711.1937}. \\[0.1 cm]
 Soit $\mathbb{K}^{3} $  le  $\mathbb{K}$ - espace vectoriel de dimension 3  o\`{u} $\mathbb{K}$  d\'{e}note 
    le corps des s\'eries  de Laurent formelles  $ \mathbb{F}_{2}((T^{-1})). $ Nous calculons en particulier des sommes exponentielles 
    (dans $\mathbb{K}^{3}$) de la forme \\
     $  \sum_{\deg Y \leq  k-1}\sum_{\deg Z \leq s-1}E(tYZ)\sum_{\deg U \leq s+m-1}E(\eta YU)\sum_{\deg V \leq s+m+l-1}E(\xi YV) $
   o\`{u} $ (t,\eta,\xi ) $ est dans la boule  unit\'{e} de $\mathbb{K}^{3}.$\\ 
    Nous d\'{e}montrons qu'elles d\'{e}pendent uniquement
   du rang de matrices triples persym\'{e}triques  avec des entr\'{e}es dans $\mathbb{F}_{2},$
 c'est-\`{a}-dire des matrices de la forme    $\left[{A\over{B \over C}}\right] $ o\`u  A  est une matrice  $ s \times k $
    persym\'{e}trique , B une matrice  $ (s+m) \times k $  persym\'{e}trique  et  C une  matrice  $ (s+m +l) \times k $  persym\'{e}trique  (une matrice \;$ [\alpha _{i,j}]  $ est 
      persym\'{e}trique  si  $ \alpha _{i,j} = \alpha _{r,s} $ \; pour  \; i+j = r+s).
        En outre, nous \'{e}tablissons plusieurs 
      formules concernant des propri\'{e}t\'{e}s de rang de partitions de matrices triples persym\'{e}triques, ce qui nous 
      conduit \`{a} une formule r\'{e}currente du nombre   $    \Gamma_{i}^{\left[s\atop{ s+m \atop s+m+l} \right]\times k} $
  des matrices de rang i de la forme  $ \left[{A\over{B \over C}}\right] $
    Nous d\'{e}duisons de cette formule r\'{e}currente 
   que si  $  0\leq i\leq\inf(s-1, k-1),$ le nombre    $  \Gamma_{i}^{\left[s\atop{ s+m \atop s+m+l} \right]\times k}  $ d\'{e}pend uniquement de i.
    D'autre part, si   $ i\geq 2s+m+1, k\geq i,\;   \Gamma_{i}^{\left[s\atop{ s+m \atop s+m+l} \right]\times k}  $ peut \^{e}tre  calcul\'{e}
   \`{a} partir du nombre  $   \Gamma_{2s' +m' +1}^{\left[s'\atop{ s'+m' \atop s'+m'+l'} \right]\times k'}   $de matrices  de rang (2s'+m'+1) de la forme  
       $ \left[{A' \over{B' \over C'}}\right]$ o\`{u}  A' est une matrice   $ s' \times k' $ persym\'{e}trique,\; B' une matrice  $ (s'+m') \times k' $ persym\'{e}trique et C'  une matrice  $ (s'+m'+l') \times k' $
        o\`{u}  s', m',l' et k' d\' {e}pendent  de  i, s, m,l et k.
          La preuve de ce r\'{e}sultat est bas\'{e}e sur une formule
        du nombre de matrices de rang i de la forme     $\left[A\over b_{-}\right] $   o\`{u}  A est une matrice double  persym\'{e}trique  et  $ b_{-} $
 une matrice ligne avec entr\'{e}es dans $ \mathbb{F}_{2}.$ Nous montrons \'{e}galement que le nombre R de repr\'{e}sentations dans 
  $ \mathbb{F}_{2}[T] $ 
 des \'{e}quations polynomiales
   \[\left\{\begin{array}{cc}
YZ + Y_{1}Z_{1} + \ldots + Y_{q-1}Z_{q-1}= 0 \\
YU + Y_{1}U_{1} + \ldots + Y_{q-1}U_{q-1}= 0\\
YV + Y_{1}V_{1} + \ldots + Y_{q-1}V_{q-1}= 0
    \end{array}\right.\]\\ 
  associ\'{e}es aux sommes exponentielles \\
   $\sum_{degY\leq  k-1}\sum_{degZ\leq s-1}E(tYZ)\sum_{\deg U\leq s+m-1}E(\eta YU) \sum_{\deg V \leq s+m+l-1}E(\xi YV)  $
  est donn\'{e} par une int\'{e}grale sur la boule unit\'{e} de $\mathbb{K}^{3}$ et est une combinaison lin\'{e}aire de
  $ \Gamma_{i}^{\left[s\atop{ s+m \atop s+m+l} \right]\times k} $ pour $ i\geq 0. $ Nous pouvons alors calculer explicitement le nombre R.
  Notre article est,  pour des raisons de longueur,  limit\' e au cas $m \geqslant 0,\; l=0.$
 \end{abstract}
  \selectlanguage{english}
  \begin{abstract}
   
   Our work concerns a generalization of the results obtained in :
  {Exponential sums and rank of  double  persymmetric  matrices over  $\mathbf{F}_2 $  }\\
{arXiv : 0711.1937}. \\[0.1 cm]
Let $ \mathbb{K}^{3}  $ be the 3-dimensional vectorspace over $ \mathbb{K}$  where $ \mathbb{K}$ denotes the field of
 Laurent Series  $ \mathbb{F}_{2}((T^{-1})). $ We compute in particular   exponential sums, (in $\mathbb{K}^{3} $)  of the form \\
  $  \sum_{\deg Y \leq  k-1}\sum_{\deg Z \leq s-1}E(tYZ)\sum_{\deg U \leq s+m-1}E(\eta YU) \sum_{\deg V \leq s+m+l-1}E(\xi YV) $ where $ (t,\eta,\xi ) $ is in the unit interval of  $\mathbb{K}^{3}.$  We  show that they only depend  on the rank  of some associated  triple persymmetric matrices  
with entries in  $\mathbb{F}_{2}, $ that is matrices of the form  $\left[{A\over{B \over C}}\right] $ 
 where A is  a  $ s \times k $ persymmetric matrix,  B a  $ (s+m) \times k $ persymmetric matrix and C is a  $ (s+m+l) \times k $ persymmetric matrix
  (A matrix $\; [\alpha _{i,j}]  $ is  persymmetric if $ \alpha _{i,j} = \alpha _{r,s}  \; for  \; i+j = r+s ).$
   Besides, we establish several formulas concerning  rank properties of  partitions of  triple  persymmetric matrices, which leads to a  recurrent formula for the number $  \Gamma_{i}^{\left[s\atop{ s+m \atop s+m+l} \right]\times k} $  of rank i  matrices of the form 
   $\left[{A\over{B \over C}}\right] $  
    We deduce from the recurrent formula that  if  $  0\leq i\leq\inf(s-1, k-1) $ then  $  \Gamma_{i}^{\left[s\atop{ s+m \atop s+m+l} \right]\times k} $  depends only on  i.
  On the other hand,  if   $ i\geq 2s+m+1, k\geq i,\;   \Gamma_{i}^{\left[s\atop{ s+m \atop s+m+l} \right]\times k}  $ 
 can be computed from the number   $   \Gamma_{2s' +m' +1}^{\left[s'\atop{ s'+m' \atop s'+m'+l'} \right]\times k'}   $    of rank (2s'+m'+1) matrices of the form  $ \left[{A' \over{B' \over C'}}\right]$ 
 where A' is  a  $ s' \times k' $ persymmetric matrix,  B' a  $ (s'+m') \times k' $ persymmetric matrix 
 and C' a  $ (s'+m'+l') \times k' $ persymmetric matrix , where s', m',l' and k' depend on i, s, m,l and k.
  The proof of this result is based on a formula  of the number of rank i matrices of the form  $\left[A\over b_{-}\right] $ 
  where A is double persymmetric and $ b_{-} $ a one-row matrix with entries in  $ \mathbb{F}_{2}.$
   We also  prove that the number R of representations in  $ \mathbb{F}_{2}[T] $ 
 of the polynomial equations
   \[\left\{\begin{array}{cc}
YZ + Y_{1}Z_{1} + \ldots + Y_{q-1}Z_{q-1}= 0 \\
YU + Y_{1}U_{1} + \ldots + Y_{q-1}U_{q-1}= 0 \\
YV + Y_{1}V_{1} + \ldots + Y_{q-1}V_{q-1}= 0
    \end{array}\right.\]\\ 
  associated to the exponential sums \\
   $\sum_{degY\leq  k-1}\sum_{degZ\leq s-1}E(tYZ)\sum_{\deg U\leq s+m-1}E(\eta YU) \sum_{\deg U\leq s+m+l-1}E(\xi YV)  $
   is given by an integral over the unit interval of $\mathbb{K}^{3}$, and is a linear combination of the
    $ \Gamma_{i}^{\left[s\atop{ s+m \atop s+m+l} \right]\times k} $ pour $ i\geq 0. $
    We can  then compute  explicitly the number R.Our article is for reasons of length limited to the case $m \geqslant 0,\; l=0$
     \end{abstract}

\newpage

  \subsection{A recurrent formula for the number of rank i matrices of the form  $\left[{A\over{B \over C}}\right] ,$ where A, B and C are persymmetric matrices over  $ \mathbb{F}_{2},$ }
\label{subsec 3.1}

   \begin{lem}
 \label{lem 3.1} Let $ s\geq 2, \; m\geq 0,\;l\geq 0, \; k\geq 1 $ and  $ 0\leq i\leq \inf{(3s+2m+l,k)}. $ Then we have the following recurrent formula
  for the number $  \Gamma_{i}^{\left[s\atop{ s+m\atop s+m+l} \right]\times k} $ of rank i matrices 
 of the form   $\left[{A\over{B \over C}}\right] $ such that  A is a  $ s\times k $ persymmetric matrix over $ \mathbb{F}_{2},$ B a 
$(s+m)\times k$ persymmetric matrix and C  a $ (s+m+l)\times k $ persymmetric  matrix  \\
  \begin{align}
  \label{eq 3.1}
& \Gamma_{i}^{\left[s\atop{ s+m\atop s+m+l} \right]\times k}  \\
& = \big[2\cdot \Gamma_{i-1}^{\left[s -1 \atop{ s+m\atop s+m+l} \right]\times k}
 +4\cdot \Gamma_{i-1}^{\left[s\atop{ s+m-1\atop s+m+l} \right]\times k}  +8\cdot\Gamma_{i-1}^{\left[s\atop{ s+m\atop s+m+l-1} \right]\times k} \big ] 
 - \big[8\cdot \Gamma_{i-2}^{\left[s -1 \atop{ s+m-1\atop s+m+l} \right]\times k}
 +16\cdot \Gamma_{i-2}^{\left[s -1\atop{ s+m\atop s+m+l-1} \right]\times k}  +32\cdot\Gamma_{i-2}^{\left[s\atop{ s+m-1\atop s+m+l-1} \right]\times k}\big] \nonumber  \\
 & + 64\cdot \Gamma_{i-3}^{\left[s -1 \atop{ s+m-1\atop s+m+l-1} \right]\times k} + \Delta _{i}^{\left[s\atop{ s+m\atop s+m+l} \right]\times k} \nonumber \\
 & \nonumber
 \end{align}
where \\

\begin{align}
  \label{eq 3.2}
  \Delta _{i}^{\left[s\atop{ s+m\atop s+m+l} \right]\times k} 
 = \sigma _{i,i,i,i}^{\left[\alpha \atop{ \beta \atop \gamma } \right]\times k}
-7\cdot\sigma _{i-1,i-1,i-1,i-1}^{\left[\alpha \atop{ \beta \atop \gamma } \right]\times k}
+14\cdot \sigma _{i-2,i-2,i-2,i-2}^{\left[\alpha \atop{ \beta \atop \gamma } \right]\times k}
-8\cdot \sigma _{i-3,i-3,i-3,i-3}^{\left[\alpha \atop{ \beta \atop \gamma } \right]\times k}
 \end{align}
 Recall that 
$$  \sigma _{i,i,i,i}^{\left[\alpha \atop{ \beta \atop \gamma } \right]\times k}$$
 denotes  the cardinality of the following set 
 $$\begin{array}{l}\Big\{ (t,\eta,\xi  ) \in \mathbb{P}/\mathbb{P}_{k+s -1}\times
           \mathbb{P}/\mathbb{P}_{k+s+m-1}\times \mathbb{P}/\mathbb{P}_{k+s+m+l-1}
\mid r(  D^{\left[s-1\atop{ s+m-1\atop s+m+l-1} \right]\times k}(t,\eta,\xi ) ) = i \\
 r(  D^{\left[s\atop{ s+m-1\atop s+m+l-1} \right]\times k}(t,\eta,\xi ) )= i,\quad  
 r(  D^{\left[s\atop{ s+m\atop s+m+l-1} \right]\times k}(t,\eta,\xi ) ) = i,  \quad
  r(  D^{\left[s\atop{ s+m\atop s+m+l} \right]\times k}(t,\eta,\xi ) ) = i \Big\}
    \end{array}$$  
 \end{lem}

\subsection{An outline of the main results in the case $m = l = 0$}
\label{subsec 3.2}
 \begin{thm}
 \label{thm 3.2}
The number  $\Gamma_{i}^{\left[s\atop{ s\atop s} \right]\times k} $ of triple persymmetric  $3s\times k $ matrices 
 over $\mathbb{F}_{2}$
  $$   \left ( \begin{array} {cccccc}
\alpha _{1} & \alpha _{2} & \alpha _{3} &  \ldots & \alpha _{k-1}  &  \alpha _{k} \\
\alpha _{2 } & \alpha _{3} & \alpha _{4}&  \ldots  &  \alpha _{k} &  \alpha _{k+1} \\
\vdots & \vdots & \vdots    &  \vdots & \vdots  &  \vdots \\
\alpha _{s-1} & \alpha _{s} & \alpha _{s +1} & \ldots  &  \alpha _{s+k-3} &  \alpha _{s+k-2}  \\
\alpha _{s} & \alpha _{s+1} & \alpha _{s +2} & \ldots  &  \alpha _{s+k-2} &  \alpha _{s+k-1}  \\
\hline \\
\beta  _{1} & \beta  _{2} & \beta  _{3} & \ldots  &  \beta_{k-1} &  \beta _{k}  \\
\beta  _{2} & \beta  _{3} & \beta  _{4} & \ldots  &  \beta_{k} &  \beta _{k+1}  \\
\vdots & \vdots & \vdots    &  \vdots & \vdots  &  \vdots \\
\beta  _{s-1} & \beta  _{s} & \beta  _{s+1} & \ldots  &  \beta_{s+k-3} &  \beta _{s+k-2}  \\
\beta  _{s} & \beta  _{s+1} & \beta  _{s+2} & \ldots  &  \beta_{s+k-2} &  \beta _{s+k-1} \\
\hline
\gamma  _{1} & \gamma   _{2} & \gamma  _{3} & \ldots  &  \gamma _{k-1} &  \gamma  _{k}  \\
\gamma   _{2} & \gamma  _{3} & \gamma   _{4} & \ldots  &  \gamma _{k} &  \gamma  _{k+1}  \\
\vdots & \vdots & \vdots    &  \vdots & \vdots  &  \vdots \\
\gamma  _{s -1} & \gamma  _{s} & \gamma  _{s+1} & \ldots  &  \gamma _{s+k-3} &  \gamma  _{s+k-2}  \\
\gamma  _{s} & \gamma  _{s+1} & \gamma  _{s+2} & \ldots  &  \gamma _{s+k-2} &  \gamma  _{s+k-1} 
\end{array}  \right) $$

of rank i is given by 
 \begin{equation}
\label{eq 3.3}
  \begin{cases}
 1  & \text{if  } i = 0, \\
 105\cdot2^{4i-6} - 21\cdot2^{3i -5}  &  \text{if  }  1\leq i\leq s-1,\;k\geq i+1, \\
  7\cdot2^{k+s-1} -7\cdot2^{2s} +105\cdot2^{4s-6} -21\cdot2^{3s-5}   &    \text{if }  i = s, \; k\geq s+1,\; s\geq 1,  \\
   147\cdot(5\cdot2^{j-1} - 1)\cdot2^{k+s+3j-6} \\
    + 21\cdot\big[5\cdot2^{4s+4j-6} - 2^{3s+3j-5} -(155\cdot2^{j-1} - 35)\cdot2^{2s+4j-7} \big] &  \text{if  }  i = s+j,  \; 1\leq j\leq s-1,\; k\geq s+j+1,\\
     7\cdot2^{2k+2s-2}+ 21\cdot\big[35\cdot2^{k+5s-7} - 39\cdot2^{k+4s-6} \big]  \\
    + 7\cdot\big[15\cdot2^{8s-6} - 465\cdot2^{7s-8} +349\cdot2^{6s-7} \big]  &  \text{if  } i=2s, \; k\geq 2s+1, \\
   105\cdot\big(2^{2k+2s+4j -2}+ 7\cdot2^{k+5s+4j-3} - 31\cdot2^{k+4s+5j-3} \big)  \\
     + 105\cdot\big(2^{8s+4j-2} - 31\cdot2^{7s+5j-3} +93\cdot2^{6s+6j-3} \big) 
     &  \text{if   }\quad  i=2s +1+j, \; k\geq 2s+2+j, \\
     &  0 \leq j\leq s-2, \\
    2^{3k+3s-3} -7\cdot2^{2k+6s-6} +7\cdot2^{k+9s-8} - 2^{12s-9}  &  \text{if  } i=3s, \; k\geq 3s  
     \end{cases}
\end{equation}  
  \begin{equation}
\label{eq 3.4}
\Gamma_{i}^{\left[s\atop{ s\atop s} \right]\times i} 
 = \begin{cases}
  2^{3s+3i-3} -7\cdot2^{4i-6} +3\cdot2^{3i-5}  &  \text{if  }\; 1\leq i\leq s+1, \\
   2^{6s+3j-3} +7\cdot2^{2s+5j-8} -7\cdot2^{2s+4j-7} -7\cdot2^{4s+4j-6} + 3\cdot2^{3s+3j-5} &  \text{if  }\;i = s+j, \;1\leq j \leq s+1,\\
     2^{9s+3j} -7\cdot2^{8s+4j-2} +7\cdot2^{7s+5j-3} - 2^{6s+6j-3}  &  \text{if  }\;i = 2s+1+j, \;0 \leq j \leq s-1,\\
  \end{cases}
\end{equation}  
  We have for $0\leq j\leq s-2,\quad k\geq 2s+2+j$\\
 \begin{align}
 \Gamma  _{2s+1+j}^{\left[s\atop{ s\atop s} \right]\times k} & = 
 16^{3j}\cdot\Gamma  _{2(s-j)+1}^{\left[s-j\atop{ s-j\atop s-j} \right]\times (k-3j)} \label{eq 3.5}
 \end{align}
 
  We have for $0\leq j\leq s-1.$\\
 \begin{align}
 \Gamma  _{2s+1+j}^{\left[s\atop{ s\atop s} \right]\times (2s+1+j)} & = 
 16^{3j}\cdot\Gamma  _{2(s-j)+1}^{\left[s-j\atop{ s-j\atop s-j} \right]\times (2(s-j)+1)} \label{eq 3.6}
  \end{align}
 We have for $k\geq 3s. $\\
  \begin{align}
 \Gamma  _{3s}^{\left[s\atop{ s\atop s} \right]\times k} & = 
 16^{3(s-1)}\cdot\Gamma  _{3}^{\left[1\atop{ 1\atop 1} \right]\times (k-3(s-1))} \label{eq 3.7}
  \end{align}

\end{thm} 
 \begin{thm}
\label{thm 3.3} Let q be a rational integer $ \geq 1,$ then \\
 \begin{align}
    g_{k,s}(t,\eta,\xi ) =  g(t,\eta,\xi ) & = \sum_{deg Y\leq k-1}\sum_{deg Z \leq  s-1}E(tYZ)\sum_{deg U \leq s-1}E(\eta YU)\sum_{deg V \leq s-1}E(\eta YV)   =
   2^{3s+k- r(D^{\left[s\atop{ s\atop s} \right]\times k}(t,\eta,\xi)) }, \label{eq 3.8} \\
   & \nonumber \\
    \int_{\mathbb{P}\times \mathbb{P}\times \mathbb{P}} g^{q}(t,\eta,\xi  )dtd\eta d\xi  
 & =2^{(3s+k)q}\cdot 2^{-3k-3s+3}\cdot \sum_{i = 0}^{\inf(3s,k)} \Gamma_{i}^{\left[s\atop{ s\atop s} \right]\times k}  \cdot2^{- qi}. \label{eq 3.9} \\
  & \nonumber 
 \end{align}
 \end{thm}
   \begin{thm}
\label{thm 3.4}
 We denote by  $ R_{q}(k,s) $ the number of solutions \\
 $(Y_1,Z_1,U_{1},V_{1}, \ldots,Y_q,Z_q,U_{q},V_{q}) $  of the polynomial equations
   \[\left\{\begin{array}{c}
 Y_{1}Z_{1} +Y_{2}Z_{2}+ \ldots + Y_{q}Z_{q} = 0,  \\
   Y_{1}U_{1} + Y_{2}U_{2} + \ldots  + Y_{q}U_{q} = 0,\\
    Y_{1}V_{1} + Y_{2}V_{2} + \ldots  + Y_{q}V_{q} = 0,\\  
 \end{array}\right.\]
  satisfying the degree conditions \\
                   $$  degY_i \leq k-1 , \quad degZ_i \leq s-1 ,\quad degU_{i}\leq s-1,\quad degV_{i}\leq s-1 \quad for \quad 1\leq i \leq q. $$ \\                           
Then \\
\begin{align}
  R_{q}(k,s) & =  \int_{\mathbb{P}\times \mathbb{P}\times \mathbb{P}} g_{k,s}^{q}(t,\eta,\xi  )dtd\eta d\xi  
 = 2^{(3s+k)q}\cdot 2^{-3k-3s+3}\cdot \sum_{i = 0}^{\inf(3s,k)} \Gamma_{i}^{\left[s\atop{ s\atop s} \right]\times k}  \cdot2^{- qi}. \label{eq 3.10}
\end{align}
  \end{thm}

\begin{example}  s = 1, $k \geq i+1$ for $0\leq i\leq 2$
  \[ \Gamma_{i}^{\left[1\atop{ 1\atop 1} \right]\times k}=
 \begin{cases}
1  &\text{if  } i = 0 \\
7\cdot(2^{k} - 1)   & \text{if   } i = 1 \\
7\cdot(2^{k}-1)\cdot(2^{k} - 2)  & \text{if   } i = 2 \\
2^{3k} - 7\cdot2^{2k} + 7\cdot2^{k+ 1}- 2^{3} & \text{if  } i = 3,\;k\geq 3
\end{cases}
\]
\end{example}

\begin{example}  s = 2, $k \geq i+1$ for $0\leq i\leq 5$
\[ \Gamma_{i}^{\left[2\atop{ 2\atop 2} \right]\times k}=
 \begin{cases}
1  &\text{if  }  i = 0 \\
21   &\text{if  }  i = 1 \\
7\cdot2^{k+ 1} + 266  & \text{if   } i = 2 \\
147\cdot2^{k+ 1} + 1344   & \text{if   } i = 3 \\
7\cdot2^{2k+2} + 651\cdot2^{k+ 2}  - 22624 & \text{if   } i = 4 \\
105\cdot2^{2k+2} - 315\cdot2^{k+ 5} + 53760   & \text{if   } i = 5 \\
2^{3k+3} - 7\cdot2^{2k+6} + 7\cdot2^{k+ 10} - 32768  & \text{if  } i = 6,\;k\geq 6
\end{cases}
\]

\end{example}
\begin{example} s = 2, k = 6.\\

 The number  $\Gamma_{i}^{\left[2\atop{ 2\atop 2} \right]\times 6}$ of rank i matrices of the form \\
  $$  \left ( \begin{array} {cccccc}
\alpha _{1} & \alpha _{2} & \alpha _{3} & \alpha _{4} & \alpha _{5} & \alpha _{6}  \\
\alpha _{2 } & \alpha _{3} & \alpha _{4} & \alpha _{5} & \alpha _{6} & \alpha _{7} \\
\beta  _{1} & \beta  _{2} & \beta  _{3}  & \beta  _{4} & \beta  _{5} & \beta  _{6}\\
\beta  _{2} & \beta  _{3} & \beta  _{4}  & \beta  _{5} & \beta  _{6} & \beta  _{7}\\
\gamma  _{1} & \gamma  _{2} & \gamma  _{3}  & \gamma  _{4} & \gamma  _{5} & \gamma  _{6} \\
\gamma  _{2} & \gamma  _{3} & \gamma  _{4}  & \gamma  _{5} & \gamma  _{6} & \gamma  _{7}
 \end{array}  \right) $$
is equal to 
\[ \begin{cases}
1  &\text{if  }  i = 0 \\
21   &\text{if  }  i = 1 \\
1162     & \text{if   } i = 2 \\
20160          & \text{if   } i = 3 \\
258720           & \text{if   } i = 4 \\
1128960           & \text{if   } i = 5 \\
688128  & \text{if  } i = 6
\end{cases}
\]
  The number of solutions \\
 $(Y_1,Z_1,U_{1},V_{1}, \ldots,Y_q,Z_q,U_{q},V_{q}) $  of the polynomial equations
   \[\left\{\begin{array}{c}
 Y_{1}Z_{1} +Y_{2}Z_{2}+ \ldots + Y_{q}Z_{q} = 0,  \\
   Y_{1}U_{1} + Y_{2}U_{2} + \ldots  + Y_{q}U_{q} = 0,\\
    Y_{1}V_{1} + Y_{2}V_{2} + \ldots  + Y_{q}V_{q} = 0,\\  
 \end{array}\right.\]
  satisfying the degree conditions \\
                   $$  degY_i \leq 5 , \quad degZ_i \leq 1 ,\quad degU_{i}\leq 1,\quad degV_{i}\leq 1 \quad for \quad 1\leq i \leq q. $$ \\                           
is equal to 
\begin{align*}
&  R_{q}(6,2)  =  \int_{\mathbb{P}\times \mathbb{P}\times \mathbb{P}} g_{6,2}^{q}(t,\eta,\xi  )dtd\eta d\xi  
 = 2^{12q-21}\cdot \sum_{i = 0}^{6} \Gamma_{i}^{\left[2\atop{ 2\atop 2} \right]\times 6}  \cdot2^{- qi}\\
& = 2^{12q-21}\cdot\big(1 + 21\cdot2^{-q} + 1162\cdot2^{-2q} + 20160\cdot2^{-3q} + 258720\cdot2^{-4q} + 1128960 \cdot2^{-5q}
 + 688128 \cdot2^{-6q}\big ) \\
 & =   2^{6q-21}\cdot\big(2^{6q} + 21\cdot2^{5q} + 1162\cdot2^{4q} + 20160\cdot2^{3q} + 258720\cdot2^{2q} + 1128960 \cdot2^{q}
 + 688128 \big ) 
\end{align*}

\end{example}

\begin{example}  s = 3, $k \geq i+1$ for $0\leq i\leq 8$
\[ \Gamma_{i}^{\left[3\atop{ 3\atop 3} \right]\times k}=
 \begin{cases}
1 &\text{if  }  i = 0 \\
21 &\text{if  }  i = 1 \\
378  &\text{if  }  i = 2 \\
7\cdot2^{k+ 2} +  5936  & \text{if   } i = 3 \\
147\cdot2^{k+ 2} + 84672  & \text{if   } i = 4 \\
147\cdot9\cdot2^{k+ 3} + 959616  & \text{if   } i = 5 \\
7\cdot2^{2k+4} + 2121\cdot2^{k+ 6} + 5863424   & \text{if   } i = 6\\
105\cdot2^{2k+4} + 2625\cdot2^{k+ 9} - 92897280  & \text{if   } i = 7 \\
105\cdot2^{2k+8} - 315\cdot2^{k+ 14} + 220200960  & \text{if   } i = 8 \\
2^{3k+6} - 7\cdot2^{2k+12} + 7\cdot2^{k+ 19} - 134217728  & \text{if  } i = 9,\;k\geq 9
\end{cases}
\]
\end{example}

 \begin{example}  s = 3, k = 5, q = 3. \\
 
 The number  $\Gamma_{i}^{\left[3\atop{ 3\atop 3} \right]\times 5}$ of rank i matrices of the form \\
  $$  \left ( \begin{array} {ccccc}
\alpha _{1} & \alpha _{2} & \alpha _{3} & \alpha _{4} & \alpha _{5}  \\
\alpha _{2 } & \alpha _{3} & \alpha _{4} & \alpha _{5} & \alpha _{6} \\
\alpha _{3} & \alpha _{4} & \alpha _{5} & \alpha _{6} & \alpha _{7} \\
\beta  _{1} & \beta  _{2} & \beta  _{3}  & \beta  _{4} & \beta  _{5}\\
\beta  _{2} & \beta  _{3} & \beta  _{4}  & \beta  _{5} & \beta  _{6}\\
\beta  _{3} & \beta  _{4} & \beta  _{5}  & \beta  _{6} & \beta  _{7}\\
\gamma  _{1} & \gamma  _{2} & \gamma  _{3}  & \gamma  _{4} & \gamma  _{5}\\
\gamma  _{2} & \gamma  _{3} & \gamma  _{4}  & \gamma  _{5} & \gamma  _{6}\\
\gamma  _{3} & \gamma  _{4} & \gamma  _{5}  & \gamma  _{6} & \gamma  _{7}
 \end{array}  \right) $$
is equal to 
 \[  \begin{cases}
1 &\text{if  }  i = 0 \\
21 &\text{if  }  i = 1 \\
378  &\text{if  }  i = 2 \\
6832  & \text{if   } i = 3 \\
103488 & \text{if   } i = 4 \\
1986432 & \text{if   } i = 5 
\end{cases}
\]

 The number of solutions \\
 $(Y_1,Z_1,U_{1},V_{1},Y_2,Z_2,U_{2},V_{2} ,Y_3,Z_3,U_{3},V_{3}) $  of the polynomial equations
   \[\left\{\begin{array}{c}
 Y_{1}Z_{1} +Y_{2}Z_{2} + Y_{3}Z_{3} = 0,  \\
   Y_{1}U_{1} + Y_{2}U_{2}  + Y_{3}U_{3} = 0,\\
    Y_{1}V_{1} + Y_{2}V_{2} + Y_{3}V_{3} = 0,\\  
 \end{array}\right.\]
  satisfying the degree conditions \\
                   $$  degY_i \leq 4 , \quad degZ_i \leq 2 ,\quad degU_{i}\leq 2,\quad degV_{i}\leq 2 \quad for \quad 1\leq i \leq 3. $$ \\                           

is equal to 
\begin{align*}
&  R_{3}(5,3)  =  \int_{\mathbb{P}\times \mathbb{P}\times \mathbb{P}} g_{5,3}^{3}(t,\eta,\xi  )dtd\eta d\xi  
 = 2^{33}\cdot \sum_{i = 0}^{5} \Gamma_{i}^{\left[3\atop{ 3\atop 3} \right]\times 5}  \cdot2^{- 3i}\\
& = 2^{33}\cdot\big(1 + 21\cdot2^{-3} + 378\cdot2^{-6} + 6832\cdot2^{-9} + 103488\cdot2^{-12} + 1986432 \cdot2^{-15}\big )  = 3563904\times 2^{18}
\end{align*}

\end{example}

\newpage
\subsection{An outline of the main results in the case $m = 1,\; l = 0$}
\label{subsec 3.3}

  \begin{thm}
 \label{thm 3.5}
   We have \\
 \begin{equation}
\label{eq 3.11}
\Gamma_{i}^{\left[s\atop{ s+1\atop s+1} \right]\times k} \\
 = \begin{cases}
 1  & \text{if  } i = 0, \\
 105\cdot2^{4i-6} - 21\cdot2^{3i -5}  &  \text{if  }  1\leq i\leq s-1,\;k\geq i+1, \\
  2^{k+s-1} -2^{2s} +105\cdot2^{4s-6} -21\cdot2^{3s-5}   &    \text{if }  i = s, \; k\geq s+1,\; s\geq 1,  \\
   33\cdot2^{k+s-1} + [105\cdot2^{4s-2} -21\cdot2^{3s-2}-69\cdot2^{2s} ]  &    \text{if }  i = s+1, \; k\geq s+2,  \\
 630\cdot2^{k+s-1} + 21\cdot[5\cdot2^{4s+2} -2^{3s+1}-65\cdot2^{2s+1} ]  &    \text{if }  i = s+2, \; k\geq s+3,  \\
  1365\cdot2^{k+s+2} + 21\cdot[5\cdot2^{4s+6} -2^{3s+4}-285\cdot2^{2s+4} ]  &    \text{if }  i = s+3, \; k\geq s+4,  \\
  2835\cdot2^{k+s+5} + 21\cdot[5\cdot2^{4s+10} -2^{3s+7}-595\cdot2^{2s+8} ]  &    \text{if }  i = s+4, \; k\geq s+5,  \\
   105\cdot(7\cdot2^{j-2} - 1)\cdot2^{k+s+3j-7} \\
    + 21\cdot\big[5\cdot2^{4s+4j-6} - 2^{3s+3j-5} - 155\cdot2^{2s+5j-10} +25\cdot2^{2s+4j-8} \big]
     &  \text{if  }  i = s+j, \\
     & \; 2\leq j\leq s,\; k\geq s+j+1,\\
    3\cdot2^{2s-1}\cdot(2^{2k}-2^{4s+4}) 
+ \big(735\cdot2^{5s-5}-393\cdot2^{4s-4}\big)\cdot\big(2^{k}-2^{2s+2}\big)  \\
+   21\cdot[5\cdot2^{8s-2} +2^{6s-4}-15\cdot2^{7s-5} ]  &    \text{if }  i = 2s+1, \; k\geq 2s+2,  \\
    53\cdot2^{2s-1}\cdot(2^{2k}-2^{4s+6}) 
+ \big(735\cdot2^{5s-1}-1629\cdot2^{4s-1}\big)\cdot\big(2^{k}-2^{2s+3}\big)  \\
+   21\cdot[5\cdot2^{8s+2} +3\cdot2^{6s}-15\cdot2^{7s} ]  &    \text{if }  i = 2s+2, \; k\geq 2s+3,  \\
 105\cdot\big(2^{2k+2s+4j +2}+ 7\cdot2^{k+5s+4j+3} - 31\cdot2^{k+4s+5j+3} \big)  \\
     + 105\cdot\big(2^{8s+4j+6} - 31\cdot2^{7s+5j+5} +93\cdot2^{6s+6j+5} \big) 
     &  \text{if   }\quad  i=2s +3+j, \; k\geq 2s+4+j, \\
     &  0 \leq j\leq s-2, \\
    2^{3k+3s-1} -7\cdot2^{2k+6s-2} +7\cdot2^{k+9s-2} - 2^{12s-1}  &  \text{if  } i=3s+2, \; k\geq 3s +2 
     \end{cases}
\end{equation}  
  \begin{equation}
\label{eq 3.12}
\Gamma_{i}^{\left[s\atop{ s+1\atop s+1} \right]\times i} 
 = \begin{cases}
  2^{3s+3i-1} -7\cdot2^{4i-6} +3\cdot2^{3i-5}  &  \text{if  }\; 1\leq i\leq s+1, \\
   2^{6s+3j-1} -7\cdot2^{4s+4j-6} +3\cdot2^{3s+3j-5} \\
    + 7\cdot2^{2s+5j-10}- 5\cdot2^{2s+4j-8}&  \text{if  }\;i = s+j, \;2\leq j \leq s+3, \\
   2^{9s+2} +7\cdot2^{7s-5}-7\cdot2^{8s-2}+7\cdot2^{6s-4} &  \text{if  }\;i = 2s+1,\\
    2^{9s+5} +7\cdot2^{7s}-7\cdot2^{8s+2}+2^{6s} &  \text{if  }\;i = 2s+2,\\
      2^{9s+8} +7\cdot2^{7s+5}-7\cdot2^{8s+6}-2^{6s+5} &  \text{if  }\;i = 2s+3,\\
     2^{9s+3j+8} -7\cdot2^{8s+4j+6} +7\cdot2^{7s+5j+5} - 2^{6s+6j+5}  &  \text{if  }\;i = 2s+3+j, \;0 \leq j \leq s-1
  \end{cases}
\end{equation} 
We have the following reduction formulas\\
 
\begin{align}
\label{eq 3.13 }
\Gamma_{2s+2+j}^{\left[s\atop{ s+1\atop s+1} \right]\times k} &
 = 16^{2j}\Gamma_{2s+2-j}^{\left[s\atop{ s+(1-j)\atop s+(1-j)} \right]\times (k-2j)} & \text{if  }  0\leq j\leq 1 \\
\Gamma_{2s+3+j}^{\left[s\atop{ s+1\atop s+1} \right]\times k} &
 = 16^{2+3j}\Gamma_{2(s-j)+1}^{\left[s -j\atop{ s -j\atop s -j} \right]\times (k-2-3j)} & \text{if  }  0\leq j\leq s-1
 \label{eq 3.14}
\end{align}
\begin{example}
We have for $s=3:$\\
 \begin{equation*}
\Gamma_{i}^{\left[3\atop{ 3+1\atop 3+1} \right]\times k} 
 = \begin{cases}
 1  & \text{if  } i = 0,\; k\geq 1, \\
 21  &  \text{if  }  i=1,\;  k\geq 2,\\
 378  &  \text{if  }  i=2, \; k\geq 3,\\ 
 2^{k+2} + 6320 &    \text{if }  i = 3, \; k\geq 4, \\
  33\cdot 2^{k+2} + 100416 &    \text{if }  i = 4, \; k\geq 5, \\ 
  630\cdot2^{k+2} + 1524096 &    \text{if }  i = 5,\;  k\geq 6, \\
   1365\cdot2^{k+5} + 21224448  &    \text{if }  i = 6, \; k\geq 7, \\
   96\cdot2^{2k} +163008\cdot2^{k+2}+1029\cdot2^{18} &    \text{if }  i = 7,\; k\geq 8,  \\
   1696\cdot2^{2k} +2176512\cdot2^{k+2}+5723\cdot2^{18} &    \text{if }  i = 8, \; k\geq 9, \\
   105\cdot2^{2k+8} +2625\cdot2^{k+15}-90720\cdot2^{18} &    \text{if }  i = 9,\; k\geq 10,  \\
    105\cdot2^{2k+12} -315\cdot2^{k+20}+215040\cdot2^{18} &    \text{if }  i = 10,\;  k\geq 11, \\
    2^{3k+8} -7\cdot2^{2k+16} +7\cdot2^{k+25} - 2^{35}  &  \text{if  } i=11. \; k\geq 11.
     \end{cases}
\end{equation*}

\end{example}
\end{thm}

\newpage
\subsection{An outline of the main results in the case $m \geqslant 2,\; l = 0$}
\label{subsec 3.4}
  \begin{thm}
 \label{thm 3.6}
 The number  $\Gamma_{i}^{\left[s\atop{ s+m\atop s+m} \right]\times k} $ of triple persymmetric  $(3s+2m)\times k $ matrices 
 over $\mathbb{F}_{2}$ of the form
 $$   \left ( \begin{array} {cccccc}
\alpha _{1} & \alpha _{2}  &  \ldots & \alpha _{k-1}  &  \alpha _{k} \\
\alpha _{2 } & \alpha _{3} &  \ldots  &  \alpha _{k} &  \alpha _{k+1} \\
\vdots & \vdots & \vdots    & \vdots  &  \vdots \\
\alpha _{s-1} & \alpha _{s} & \ldots  &  \alpha _{s+k-3} &  \alpha _{s+k-2}  \\
 \alpha  _{s } & \alpha  _{s +1} & \ldots & \alpha  _{s +k-2}& \alpha  _{s +k-1}\\
\hline \\
\beta  _{1} & \beta  _{2}  & \ldots  &  \beta_{k-1} &  \beta _{k}  \\
\beta  _{2} & \beta  _{3}  & \ldots  &  \beta_{k} &  \beta _{k+1}  \\
\vdots & \vdots    &  \vdots & \vdots  &  \vdots \\
\beta  _{m+1} & \beta  _{m+2}  & \ldots  &  \beta_{k+m-1} &  \beta _{k+m}  \\
\vdots & \vdots    &  \vdots & \vdots  &  \vdots \\
\beta  _{s+m-1} & \beta  _{s+m}  & \ldots  &  \beta_{s+m+k-3} &  \beta _{s+m+k-2}  \\
 \beta _{s+m} & \beta _{s+m+1} & \ldots & \beta _{s+m+k-2} & \beta _{s+m+k-1}\\
\hline \\
\gamma  _{1} & \gamma   _{2}  & \ldots  & \gamma  _{k-1} &  \gamma  _{k}  \\
\gamma  _{2} & \gamma  _{3}  & \ldots  & \gamma  _{k} &  \gamma  _{k+1}  \\
\vdots & \vdots    &  \vdots & \vdots  &  \vdots \\
 \gamma  _{m+1} &  \gamma _{m+2}  & \ldots  & \gamma _{k+m-1} &  \gamma  _{k+m}  \\
\vdots & \vdots   &  \vdots & \vdots  &  \vdots \\
 \gamma  _{s+m-1} & \gamma  _{s+m}  & \ldots  & \gamma  _{s+m+k-3} &  \gamma  _{s+m+k-2}  \\
  \gamma  _{s+m} & \gamma  _{s+m+1}  & \ldots  & \gamma  _{s+m+k-2} &  \gamma  _{s+m+k-1}  
\end{array}  \right). $$ is given by

  \begin{equation}
    \label{eq 3.15}  
    \begin{cases}
 1  & \text{if  } i = 0,\;k\geq 1 \\
 105\cdot2^{4i-6} - 21\cdot2^{3i -5}  &  \text{if  }  1\leq i\leq s-1,\;k\geq i+1, \\ 
2^{k+s-1} -2^{2s}+21\cdot\big(5\cdot2^{4s-6}-2^{3s-5}\big)  &\text{if  }  i = s,\;k\geq s+1\\
& \\
(21\cdot2^{j-1}-3)\cdot2^{k+s+2j-4}\\
 +21\cdot\big(5\cdot2^{4s+4j-6}-2^{3s+3j-5}-5\cdot2^{2s+4j-6}+2^{2s+3j-5}\big) 
  &\text{if  } i=s+j,\;  1\leq j\leq m-1 ,\\
  & k\geq s+j+1\\
  & \\
   \big(21\cdot2^{s+3m-5}+45\cdot2^{s+2m-4}\big)\cdot\big(2^{k}-2^{s+m+1}\big)\\
 + 105\cdot2^{4s+4m-6}-21\cdot2^{3s+3m-5}-21\cdot2^{2s+4m-6}+9\cdot2^{2s+3m-5} 
   &\text{if  }  i=s+m ,\;k\geq s+m+1\\
   & \\
    21\cdot\big[ 2^{k+s+3m+3j-5}+35\cdot2^{k+s+2m+4j-7}-9\cdot2^{k+s+2m+3j-6}\\
       + 5\cdot2^{4s+4m+4j-6} -2^{3s+3m+3j-5}- 5\cdot2^{2s+4m+4j-6}\\
       -155\cdot2^{2s+3m+5j-8}+45\cdot2^{2s+3m+4j-7}\big ]
   &\text{if  } i=s+m+j,\;1\leq j\leq s-1,\\
   & k\geq s+m+j+1\\
   & \\
   3\cdot2^{2k+2s+m-2}+21\cdot2^{k+4s+3m-5}+735\cdot2^{k+5s+2m-7} - 477\cdot2^{k+4s+2m-6}\\
 +105\cdot2^{8s+4m-6} -105\cdot2^{6s+4m-6}-3255\cdot2^{7s+3m-8}+1629\cdot2^{6s+3m-7} 
  &\text{if  } i=2s+m ,\;k\geq 2s+m+1\\ 
  & \\
    21\cdot2^{2k+2s+m+3j-2}+21\cdot2^{k+4s+3m+3j-2}+735\cdot2^{k+5s+2m+4j-3}\\
      - 945\cdot2^{k+4s+2m+4j-3} +105\cdot2^{8s+4m+4j-2} -105\cdot2^{6s+4m+4j-2}\\
      -3255\cdot2^{7s+3m+5j-3}+3255\cdot2^{6s+3m+5j-3}
   &\text{if  }i=2s+m+1+j,\\
   & 0\leq j\leq m-2, \;  k\geq 2s+m+2+j \\
   & \\
    53\cdot2^{2s-1}\cdot(2^{2k+4m-4}-2^{4s+8m-2}) \\
 + \big(735\cdot2^{5s-1}-1629\cdot2^{4s-1}\big)\cdot\big(2^{k+6m-6}-2^{2s+8m-5}\big)  \\
+   21\cdot[5\cdot2^{8s+8m-6} +3\cdot2^{6s+8m-8}-15\cdot2^{7s+8m-8} ]  &  \text{if } i=2s+2m,\;
  k\geq 2s+2m+1 \\
 & \\
  105\cdot\big(2^{2k+2s+4m+4j-2}+ 7\cdot2^{k+5s+6m+4j-3} - 31\cdot2^{k+4s+6m+5j-3} \big)  \\
    + 105\cdot\big(2^{8s+8m+4j-2} - 31\cdot2^{7s+8m+5j-3} +93\cdot2^{6s+8m+6j-3} \big)  
 &\text{if  }i=2s+2m+1+j, \\
 &  0\leq j\leq s-2,\; k\geq 2s+2m+2+j\\
 & \\
  2^{3k+2m+3s-3}-7\cdot2^{2k+4m+6s-6}+7\cdot2^{k+6m+9s-8}-2^{8m+12s-9} &\text{if  } i=3s+2m,\; k\geq 3s+2m
 \end{cases}
\end{equation}
 \newpage
 \begin{equation}
 \label{eq 3.16}
   \Gamma_{i}^{\left[s\atop{ s +m\atop s+m} \right]\times i} 
 \end{equation}
 \begin{equation*}  
  =  \begin{cases}
 2^{3s+2m+3i-3}-7\cdot2^{4i-6}+3\cdot2^{3i-5}  &\text{if  }  1\leq i\leq s+1\\
 & \\
  2^{6s+2m+3j-3}-7\cdot2^{4s+4j-6}+3\cdot2^{3s+3j-5}\\
  +(2^{j-1}-1)\cdot2^{2s+3j-5} &\text{if  }i=s+j,\; 1\leq j\leq m+1  \\
  & \\
   2^{6s+5m+3j-3}-7\cdot2^{4s+4m+4j-6}+3\cdot2^{3s+3m+3j-5}\\
+ 2^{2s+4m+4j-6}+7\cdot2^{2s+3m+5j-8}- 9\cdot2^{2s+3m+4j-7}&\text{if  } i=s+m+j,\;1\leq j\leq s-1\\
& \\
 2^{9s+5m-3}-7\cdot2^{8s+4m-6}+7\cdot2^{7s+3m-8}+3\cdot2^{6s+3m-7}+2^{6s+4m-6} &\text{if  } i=2s+m\\
  & \\
 2^{9s+5m+2j}-7\cdot2^{8s+4m+4j-2}+7\cdot2^{7s+3m+5j-3}-3\cdot2^{6s+3m+5j-3}+2^{6s+4m+4j-2} &\text{if  }i=2s+m+1+j,\\
 &  0\leq j\leq m-2\\
 & \\
  2^{9s+8m-3}-7\cdot2^{8s+8m-6}+7\cdot2^{7s+8m-8}-2^{6s+8m-8} &\text{if  } i=2s+2m\\ 
  & \\
 2^{9s+8m+3j}-7\cdot2^{8s+8m+4j-2}+7\cdot2^{7s+8m+5j-3}-2^{6s+8m+6j-3} &\text{if  }i=2s+2m+1+j,\\
 &  0\leq j\leq s-2 \\
 & \\
  21\cdot2^{8m+12s-9} &\text{if  } i=3s+2m
 \end{cases}
\end{equation*}
  We have for $ 0\leq j\leq m-2,\;k \geq 2s+m+2+j$\\[0.1 cm]
  
  \begin{align}
 \Gamma_{2s+m+1+j}^{\left[s\atop{ s+m\atop s+m} \right]\times k}
 = 16^{2j}\Gamma_{2s+1+(m-j)}^{\left[s\atop{ s+(m-j)\atop s+(m-j)} \right]\times (k-2j)} \label{eq 3.17}\\
 & \nonumber
       \end{align}
   We have for $ j=m-1,\;k \geq 2s+2m+1$\\[0.1 cm]
   
    \begin{align}
 \Gamma_{2s+2m}^{\left[s\atop{ s+m\atop s+m} \right]\times k}
 = 16^{2m-2}\Gamma_{2s+2}^{\left[s\atop{ s+1\atop s+1} \right]\times (k-2(m-1))} \label{eq 3.18}\\
 & \nonumber
 \end{align}

 We have for $ 0\leq j\leq s-2,\;k \geq 2s+2m+2+j$\\[0.1 cm]
  \begin{align}
  \Gamma_{2s+2m+1+j}^{\left[s\atop{ s+m\atop s+m} \right]\times k} 
 = 16^{2m+3j}\Gamma_{2(s-j)+1}^{\left[s -j\atop{ s -j\atop s -j} \right]\times (k-2m-3j)} \label{eq 3.19} \\
 & \nonumber
  \end{align}
  
  \newpage
We have for $j=s-1,\;k\geq 3s+2m$\\[0.1 cm]

   \begin{align}
  \Gamma_{3s+2m}^{\left[s\atop{ s+m\atop s+m} \right]\times k} 
 = 16^{2m+3s-3}\Gamma_{3}^{\left[1\atop{ 1\atop 1} \right]\times (k-2m-3s+3)}  \label{eq 3.20}
   \end{align}

 \begin{align}
 & \text{ Let q be a rational integer $ \geq 1,$ then }\label{eq 3.21} \\
  &  g_{k,s,m}(t,\eta,\xi ) =  g(t,\eta,\xi )  = \sum_{deg Y\leq k-1}\sum_{deg Z \leq  s-1}E(tYZ)\sum_{deg U \leq s+m-1}E(\eta YU)\sum_{deg V \leq s+m-1}E(\eta YV) \nonumber \\
 & =  2^{3s+2m+k- r(D^{\left[s\atop{ s+m\atop s+m} \right]\times k}(t,\eta,\xi)) }\nonumber \\
 & \nonumber \\
 &   \quad   \text{and} \nonumber \\
  & \nonumber \\
 &\int_{\mathbb{P}^{3}} g^{q}(t,\eta,\xi  )d t d \eta d\xi = \nonumber\\
& = \sum_{(t,\eta,\xi  )\in \mathbb{P}/\mathbb{P}_{k+s-1}\times \mathbb{P}/\mathbb{P}_{k+s+m-1}\times \mathbb{P}/\mathbb{P}_{k+s+m-1}\atop { r( D^{\left[\stackrel{s}{\stackrel{s+m}{s+m}}\right] \times k }(t,\eta ,\xi )  =  i}}
 2^{\big(k+3s+2m-  r( D^{\left[\stackrel{s}{\stackrel{s+m}{s+m}}\right] \times k }(t,\eta ,\xi ) )\big)q  } \int_{\mathbb{P}_{k+s-1}}dt \int_{\mathbb{P}_{k+s+m-1}}d\eta \int_{\mathbb{P}_{k+s+m-1}}d\xi \nonumber \\
 & = \sum_{i = 0}^{\inf(k,3s+2m)}\sum_{(t,\eta,\xi  )\in \mathbb{P}/\mathbb{P}_{k+s-1}\times \mathbb{P}/\mathbb{P}_{k+s+m-1}\times \mathbb{P}/\mathbb{P}_{k+s+m-1}}
 2^{(k+3s+2m-i)q  } \int_{\mathbb{P}_{k+s-1}}dt \int_{\mathbb{P}_{k+s+m-1}}d\eta \int_{\mathbb{P}_{k+s+m-1}}d\xi \nonumber \\
  & = 2^{(k+3s+2m)q -(3k +3s+2m-3)}
 \sum_{i=0}^{\inf(k,3s+2m)}\Gamma_{i}^{\left[s\atop{ s+m \atop s+m }\right]\times k}\cdot2^{-iq} \nonumber
 & \nonumber \\
 & \nonumber
 \end{align}  

 We denote by  $ R_{q}(k,s,m) $ the number of solutions \\
 $(Y_1,Z_1,U_{1},V_{1}, \ldots,Y_q,Z_q,U_{q},V_{q}) $  of the polynomial equations
   \[\left\{\begin{array}{c}
 Y_{1}Z_{1} +Y_{2}Z_{2}+ \ldots + Y_{q}Z_{q} = 0,  \\
   Y_{1}U_{1} + Y_{2}U_{2} + \ldots  + Y_{q}U_{q} = 0,\\
    Y_{1}V_{1} + Y_{2}V_{2} + \ldots  + Y_{q}V_{q} = 0,\\  
 \end{array}\right.\]
  satisfying the degree conditions \\
                   $$  degY_i \leq k-1 , \quad degZ_i \leq s-1 ,\quad degU_{i}\leq s+m-1,\quad degV_{i}\leq s+m-1 \quad for \quad 1\leq i \leq q. $$ \\                           
Then \\
\begin{align}
  R_{q}(k,s,m) & =  \int_{\mathbb{P}\times \mathbb{P}\times \mathbb{P}} g_{k,s,m}^{q}(t,\eta,\xi  )dtd\eta d\xi  
  = 2^{(k+3s+2m)q -(3k +3s+2m-3)}
 \sum_{i=0}^{\inf(k,3s+2m)}\Gamma_{i}^{\left[s\atop{ s+m \atop s+m }\right]\times k}\cdot2^{-iq} \label{eq 3.20}
\end{align}

  \begin{example}
We have for $s=3,\;m=4,\;k=10:$\\
 \begin{equation*}
\Gamma_{i}^{\left[3\atop{ 3+4\atop 3+4} \right]\times 10} 
 = \begin{cases}
 1  & \text{if  } i = 0 \\
 21  &  \text{if  }  i=1\\
 378  &  \text{if  }  i=2\\ 
 10416 &    \text{if }  i = 3\\
  140352 &    \text{if }  i = 4,  \\ 
  1994112 &    \text{if }  i = 5, \\
   29598720 &    \text{if }  i = 6 \\
   458661888 &    \text{if }  i = 7, \\
   109389\cdot2^{16} &    \text{if }  i = 8, \\
   213759\cdot2^{19} &    \text{if }  i = 9, \\
  2^{44} -14273\cdot2^{23}  &  \text{if  } i=10
     \end{cases}
\end{equation*}

\end{example}

\begin{example} s = 3, m=4, k = 7, q = 3\\

 The number  $\Gamma_{i}^{\left[3\atop{ 3+4\atop 3+4} \right]\times 7}$ of rank i matrices of the form \\
  $$  \left ( \begin{array} {ccccccc}
\alpha _{1} & \alpha _{2} & \alpha _{3} & \alpha _{4} & \alpha _{5} & \alpha _{6} & \alpha _{7} \\
\alpha _{2 } & \alpha _{3} & \alpha _{4} & \alpha _{5} & \alpha _{6} & \alpha _{7}& \alpha _{8} \\
\alpha _{3 } & \alpha _{4} & \alpha _{5} & \alpha _{6} & \alpha _{7} & \alpha _{8}& \alpha _{9} \\
\beta  _{1} & \beta  _{2} & \beta  _{3}  & \beta  _{4} & \beta  _{5} & \beta  _{6} & \beta  _{7}\\
\beta  _{2} & \beta  _{3} & \beta  _{4}  & \beta  _{5} & \beta  _{6} & \beta  _{7} & \beta  _{8}\\
\beta  _{3} & \beta  _{4} & \beta  _{5}  & \beta  _{6} & \beta  _{7} & \beta  _{8} & \beta  _{9}\\
\beta  _{4} & \beta  _{5} & \beta  _{6}  & \beta  _{7} & \beta  _{8} & \beta  _{9} & \beta  _{10}\\
\beta  _{5} & \beta  _{6} & \beta  _{7}  & \beta  _{8} & \beta  _{9} & \beta  _{10} & \beta  _{11}\\
\beta  _{6} & \beta  _{7} & \beta  _{8}  & \beta  _{9} & \beta  _{10} & \beta  _{11} & \beta  _{12}\\
\beta  _{7} & \beta  _{8} & \beta  _{9}  & \beta  _{10} & \beta  _{11} & \beta  _{12} & \beta  _{13}\\
\gamma  _{1} & \gamma  _{2} & \gamma  _{3}  & \gamma  _{4} & \gamma  _{5} & \gamma  _{6}  & \gamma  _{7}\\
\gamma  _{2} & \gamma  _{3} & \gamma  _{4}  & \gamma  _{5} & \gamma  _{6} & \gamma  _{7} & \gamma  _{8}\\
\gamma  _{3} & \gamma  _{4} & \gamma  _{5}  & \gamma  _{6} & \gamma  _{7} & \gamma  _{8}  & \gamma  _{9}\\
\gamma  _{4} & \gamma  _{5} & \gamma  _{6}  & \gamma  _{7} & \gamma  _{8} & \gamma  _{9} & \gamma  _{10}\\
\gamma  _{5} & \gamma  _{6} & \gamma  _{7}  & \gamma  _{8} & \gamma  _{9} & \gamma  _{10}  & \gamma  _{11}\\
\gamma  _{6} & \gamma  _{7} & \gamma  _{8}  & \gamma  _{9} & \gamma  _{10} & \gamma  _{11} & \gamma  _{12}\\
\gamma  _{7} & \gamma  _{8} & \gamma  _{9}  & \gamma  _{10} & \gamma  _{11} & \gamma  _{12} & \gamma  _{13}
 \end{array}  \right) $$
is equal to

\[ \begin{cases}
 1  & \text{if  } i = 0 \\
 21  &  \text{if  }  i=1\\
 378  &  \text{if  }  i=2\\ 
 6832 &    \text{if }  i = 3\\
  108096 &    \text{if }  i = 4,  \\ 
  1714560 &    \text{if }  i = 5, \\
   27276288 &    \text{if }  i = 6 \\
  2^{35} -3553\cdot2^{13} &    \text{if }  i = 7, \\
      \end{cases}\]

 The number of solutions \\
 $(Y_1,Z_1,U_{1},V_{1},Y_2,Z_2,U_{2},V_{2} ,Y_3,Z_3,U_{3},V_{3}) $  of the polynomial equations
   \[\left\{\begin{array}{c}
 Y_{1}Z_{1} +Y_{2}Z_{2} + Y_{3}Z_{3} = 0,  \\
   Y_{1}U_{1} + Y_{2}U_{2}  + Y_{3}U_{3} = 0,\\
    Y_{1}V_{1} + Y_{2}V_{2} + Y_{3}V_{3} = 0,\\  
 \end{array}\right.\]
  satisfying the degree conditions \\
                   $$  degY_i \leq 6 , \quad degZ_i \leq 2 ,\quad degU_{i}\leq 6,\quad degV_{i}\leq 6 \quad for \quad 1\leq i \leq 3. $$ \\                           

is equal to 
\begin{align*}
&  R_{3}(7,3,4)  =  \int_{\mathbb{P}\times \mathbb{P}\times \mathbb{P}} g_{7,3,4}^{3}(t,\eta,\xi  )dtd\eta d\xi  
 = 2^{37}\cdot \sum_{i = 0}^{7} \Gamma_{i}^{\left[3\atop{ 3+4\atop 3+4} \right]\times 7}  \cdot2^{- 3i}\\
 & = 2^{37}\cdot\big(1 + 21\cdot2^{-3} + 378\cdot2^{-6} + 6832\cdot2^{-9} + 108096 \cdot 2^{-12} + 1714560 \cdot2^{-15}\\
& +  27276288 \cdot 2^{-18}  +  (2^{35} -3553\cdot2^{13}) \cdot 2^{ -21}  \big )  = 4243395\cdot2^{29}
\end{align*}
\end{example}

 \begin{example} The fraction of square triple persymmetric   $ \left[s\atop{ s+m\atop s+m} \right]\times (3s+2m) $   matrices which are  invertible is equal to 
   $ \displaystyle \frac{\Gamma_{3s+2m}^{\left[s\atop{ s+m\atop s+m} \right]\times (3s+2m)}}
 {\sum_{i = 0}^{3s+2m}\Gamma_{i}^{\left[s\atop{ s+m\atop s+m} \right]\times (3s+2m)} } = \frac{21}{64}. $
   \end{example}

 \end{thm}

 \end{document}